\documentclass[3p, 11pt]{elsarticle}

\usepackage{color}
\usepackage[colorlinks,linkcolor=blue,citecolor=blue]{hyperref}
\usepackage{amsmath,amsfonts,amssymb,amsthm,stmaryrd,bm}
\usepackage{mathrsfs} 
\usepackage{graphicx,epstopdf,float,subfigure}
\usepackage[margin=10pt,font=small,labelfont=bf,labelsep=period]{caption}

\usepackage{multirow,booktabs}

\numberwithin{equation}{section} 
\newtheorem{theorem}{Theorem}[section]
\newtheorem{lemma}{Lemma}[section]

\newtheorem{remark}{Remark}[section]
\newtheorem{example}{Example}[section]
\biboptions{sort&compress}

\def \d {\mathrm{d}}

\newcommand{\normm}[1]{\interleave  #1  \interleave}

\begin{document}
	
\begin{frontmatter}
	
\title{The interior penalty virtual element method for fourth-order singular perturbation problems}

\author{Fang Feng}
\ead{ffeng@njust.edu.cn}
\address{School of Mathematics and Statistics, Nanjing University of Science and Technology, Nanjing, China}
\author{Yue Yu\footnote{Corresponding author.}}
\ead{terenceyuyue@sjtu.edu.cn}
\address{School of Mathematics and Computational Science, Hunan Key Laboratory for Computation and Simulation in Science and Engineering, Key Laboratory of Intelligent Computing and Information Processing of Ministry of Education, Xiangtan University, Xiangtan, Hunan 411105, PR China}


\begin{abstract}
This paper is dedicated to the numerical solution of a fourth-order singular perturbation problem using the interior penalty virtual element method (IPVEM) proposed in \cite{ZMZW2023IPVEM}. The study introduces modifications to the jumps and averages in the penalty term, as well as presents an automated mesh-dependent selection of the penalty parameter. Drawing inspiration from the modified Morley finite element methods, we leverage the conforming interpolation technique to handle the lower part of the bilinear form. Through our analysis, we establish optimal convergence in the energy norm and provide a rigorous proof of uniform convergence concerning the perturbation parameter in the lowest-order case.
\end{abstract}
	

\end{frontmatter}

%
%
\section{Introduction}

	Let $\Omega$ be a bounded  polygonal domain of $\mathbb{R}^2$ with boundary $\partial \Omega$.  For $f\in L^2(\Omega)$, we consider the following boundary value problem of the fourth-order singular perturbation equation:
	\begin{equation}\label{originalSP}
		\left\{\begin{array}{ll}
			\varepsilon^2\Delta^2u-\Delta u = f & {\rm in}~~\Omega,\\
			u = \dfrac{\partial u}{\partial \bm n}=0& {\rm on}~~\partial \Omega,
		\end{array}\right.
	\end{equation}
	where $\bm{n} = (n_1,n_2)$ is the unit outer normal to $\partial \Omega$, $\Delta$ is the standard Laplacian operator and $\varepsilon$ is a real parameter satisfying $0< \varepsilon\le 1$.
When $\varepsilon$ is not small, the problem can be numerically resolved as the biharmonic equation by using some standard finite element methods for fourth-order problems.
However, in this article, we are primarily concerned with the case where $\varepsilon \to 0$, or the differential equation formally degenerates to the Poisson equation.
Due to the presence of the fourth-order term,  $H^2$-conforming finite elements with $C^1$ continuity were considered in \cite{Semper1992}, but these methods are quite complicated even in two dimensions. To overcome this difficulty, it is prefer to using nonconforming finite elements. Since the differential equation reduces to the Poisson equation in the singular limit, $C^0$-nonconforming finite elements are better suited for the task, as shown in \cite{NTW01}. In that work, the method was proved to converge uniformly in the perturbation parameter in the energy norm, and a counterexample was given to show that the Morley method diverges for the reduced second-order equation (see also the patch test in \cite{Wang01}).  Considering that the Morley element has the least number of element degree of freedom (DoF), Wang et al. proposed in \cite{WXH06} a modified Morley element method. This method still uses the DoFs of the Morley element, but linear approximation of finite element functions is used in the lower part of the bilinear form:
\begin{equation}\label{mMorley}
\varepsilon^2  a_h(u_h,v_h) + b_h(\Pi_h^1 u_h, \Pi_h^1 v_h) = (f, \Pi_h^1 v_h),
\end{equation}
where
\[a_h(v_h,w_h) = \sum\limits_{K \in \mathcal{T}_h} (\nabla^2 v_h, \nabla^2 w_h)_K, \qquad
b_h(v_h, w_h) = \sum\limits_{K \in \mathcal{T}_h} (\nabla v_h, \nabla w_h)_K,\]
and $\Pi_h^1$ is the interpolation operator corresponding to linear conforming element for second-order partial differential equations. It was also shown that the modified method converges uniformly in the perturbation parameter.
On the other hand, among the parameter-robust methods, the $C^0$ interior penalty methods \cite{Brenner2005SungC0IP,Brenner2011NeilanC0IP} may be most attractive as they utilize the standard $C^0$ finite element spaces. These spaces are primarily designed for solving second-order problems and effectively reduce the computational cost to a certain extent.

This paper focuses on virtual element methods (VEMs), which are a generalization of the standard finite element method that allows for general polytopal meshes. First proposed and analyzed in \cite{Beirao-Brezzi-Cangiani-2013}, with other pioneering works found in \cite{Ahmad-Alsaedi-Brezzi-2013,Beirao-Brezzi-Marini-2014}, VEMs have several advantages over standard finite element methods. For example, they are more convenient for handling partial differential equations on complex geometric domains or those associated with high-regularity admissible spaces. To date, a considerable number of conforming and nonconforming VEMs have been developed to address second-order elliptic equations \cite{Ahmad-Alsaedi-Brezzi-2013,Beirao-Brezzi-Cangiani-2013,DeDios-Lipnikov-Manzini-2016,Cangiani-Manzini-Sutton-2016} as well as fourth-order elliptic equations \cite{Brezzi-Marini-2013,Chinosi-Marini-2016,Antonietti-Manzini-Verani-2018,Zhao-Zhang-Chen-2018,Chen-HuangX-2020}.
The versatility and effectiveness of the VEM have made it a popular choice for tackling a wide range of scientific and engineering problems, including the Stokes problem \cite{Antonietti-Bruggi-Scacchi-17,Caceres-Gatica-2017,Liu-Li-Chen-2017}, the Navier-Stokes equation \cite{Beirao-Lovadina-Vacca-2018,Gatica-Munar-2018,Beirao-Mora-Vacca-2019}, the MHD equations \cite{Alvarez-Beirao-Dassi-2023,Beirao-Dassi-Manzini-2023}, the sine-Gordon equation \cite{Adak-Nataraj-2020}, the topology optimization problems \cite{Chi-Pereira-Menezes-Paulino-20,Zhang-chi-Paulino-20,Antonietti-Bruggi-Scacchi-17} and the variational and hemivariational inequalities \cite{FHH22,Wang-Wei-2017,FHH19,Wang-wu-han-2021,Ling-wang-han-2020,qiu-wang-ling-zhao-2023,xiao-ling-2023,wang-zhao-2021}.
For a comprehensive understanding of recent advancements in the VEM, we recommend referring to the book \cite{Antonietti-Beirao-Manzini2022} and the associated references. Additionally, with regards to the fourth-order singular perturbation problems, a notable application of the
$C^0$-continuous nonconforming virtual element method can be found in \cite{ZZC20}.

Recently, the interior penalty technique for VEMs has been explored in \cite{ZMZW2023IPVEM} for the biharmonic equation, which is equipped with the same DoFs for the $H^1$-conforming virtual elements, and the new numerical scheme~---~the interior penalty virtual element method (IPVEM)~---~can be regarded as a combination of the virtual element space and discontinuous Galerkin scheme, since the resulting global discrete space is not $C^0$-continuous and an interior penalty formulation is adopted to enforce the $C^1$ continuity of the solution. Inspired by the technique in \eqref{mMorley} for discontinuous elements, we are trying to explore the feasibility of the IPVEM for solving the fourth-order singular perturbation problem \eqref{originalSP}. In the context of VEMs, however, the linear conforming interpolation operator can be omitted as it yields the same elliptic projection.
	
	The remainder of the paper is structured as follows. We begin by introducing the continuous variational problem and presenting some useful results in VEM analysis in Section \ref{sec:cvariationalProb}.
Section \ref{sec:IPVEM} is dedicated to the introduction of the IPVEM. In contrast to the jump and average terms in \cite{ZMZW2023IPVEM},  we include the elliptic projector $\Pi_h^\nabla$ for all $v$ and $w$ in the penalty terms $J_i(v,w)$ for $i=1,2,3$ to simplify the implementation. Additionally, we provide an automated mesh-dependent selection of the penalty parameter $\lambda_e$, following a similar deduction as described in \cite{Carsten2023IP}.
In Section \ref{sec:err}, we establish the optimal convergence of the IPVEM in the energy norm and provide a uniform error estimate in the lowest-order case. Our analysis demonstrates that the IPVEM is robust with respect to the perturbation parameter. This is based on the observation that we can equivalently include the $H^1$-conforming interpolation in the lower part of the bilinear form, as in \eqref{mMorley}, since all required degrees of freedom are accessible.
Numerical examples are presented in Section \ref{sec:numerical} to validate the theoretical predictions. Finally, conclusions are provided in the last section.

\section{The continuous variational problem} \label{sec:cvariationalProb}

We first introduce some notations and symbols frequently used in this paper. For a bounded Lipschitz domain $D$, the symbol $( \cdot , \cdot )_D$ denotes the $L^2$-inner product on $D$, $\|\cdot\|_{0,D}$ denotes the $L^2$-norm, and $|\cdot|_{s,D}$ is the $H^s(D)$-seminorm. For all integer $k\ge 0$, $\mathbb{P}_k(D)$ is the set of polynomials of degree $\le k$ on $D$. Let $e \subset \partial K$ be the common edge for elements $K = K^-$ and $K^+$, and let $v$ be a scalar function defined on $e$. We introduce the jump and average of $v$ on $e$ by $[v] = v^- - v^+$ and $\{v\} = \frac12 (v^- + v^+)$, where $v^-$ and $v^+$ are the traces of $v$ on $e$ from the interior and exterior of $K$, respectively. On a boundary edge, $[ v ] = v$ and $\{ v \} = v $. Moreover, for any two quantities $a$ and $b$, ``$a\lesssim b$" indicates ``$a\le C b$" with the constant $C$ independent of the mesh size $h_K$, and ``$a\eqsim b$" abbreviates ``$a\lesssim b\lesssim a$".

 The variational formulation of \eqref{originalSP} reads: Find $u\in V:=H_0^2(\Omega)$ such that
\begin{equation}\label{origion}
	\varepsilon^2a(u,v)+b(u,v)=(f,v),\quad  v\in H_0^2(\Omega),
\end{equation}
where
$a(u,v) = (\nabla^2u,\nabla^2v)$ and  $ b(u,v)=(\nabla u,\nabla v) $.
To avoid complicated presentation, we confine our discussion in two dimensions, with the family of polygonal meshes  $\{ \mathcal{T}_h \}_{h>0}$ satisfying the following condition (cf. \cite{Brezzi-Buffa-Lipnikov-2009,Chen-HuangJ-2018}):
\begin{enumerate}[{\bf H}.]
\item For each $K\in {\mathcal T}_h$, there exists a ``virtual triangulation" ${\mathcal T}_K$ of $K$ such that ${\mathcal T}_K$ is uniformly shape regular and quasi-uniform. The corresponding mesh size of ${\mathcal T}_K$ is proportional to $h_K$. Each edge of $K$ is a side of a certain triangle in ${\mathcal T}_K$.
\end{enumerate}
As shown in \cite{Chen-HuangJ-2018}, this condition covers the usual geometric assumptions frequently used in the context of VEMs.
Under this geometric assumption, we can establish some fundamental results in VEM analysis as used in \cite{Huang2021YuMedius}.

According to the standard Dupont-Scott theory (cf. \cite{BS2008}), for all $v\in H^l(K)$ ($0\le l \le k$) there exists
a certain $q \in \mathbb{P}_{l-1}(K)$ such that
\begin{equation}\label{BHe1}
|v - q|_{m,K} \lesssim h_K^{l - m} | v |_{l,K},\quad  m\le l.
\end{equation}
The following trace inequalities are very useful for our forthcoming analysis.
\begin{lemma}\label{lem:trace}
For any $\varepsilon  > 0$, there exists a constant $C(\varepsilon )$ such that
\begin{align*}
&\|v\|_{0,\partial K}^2 \lesssim \|v\|_{0,K} \|v\|_{1,K}, \quad v \in H^1(K), \\
&\| v \|_{0,\partial K} \lesssim \varepsilon h_K^{1/2}| v |_{1,K} + C(\varepsilon )h_K^{ - 1/2}\| v \|_{0,K},\quad v \in H^1(K). \label{trace1}
\end{align*}
\end{lemma}

\section{The interior penalty virtual element method} \label{sec:IPVEM}

\subsection{The $H^1$-projection in the lifting space}

The interior penalty virtual element method (IPVEM) was proposed in \cite{ZMZW2023IPVEM}. In the construction, the authors first introduced a $C^1$-conforming virtual element space
\[\left\{v \in H^2(K):  \Delta^2 v \in \mathbb{P}_k(K), v|_e \in \mathbb{P}_{k+2}(e), \partial_{\bm n} v|_e \in \mathbb{P}_{k+1}(e), ~~e \subset \partial K\right\}, \quad k \ge 2,\]
which can be viewed as the lifting version of the standard $k$-th order $C^1$-conforming virtual element space for fourth-order singular perturbation problems. By checking the computability, we find that the order of the interior moments can be actually reduced to $k-2$. To this end, we instead consider the following modified lifting $C^1$-conforming virtual element space
\[
\widetilde{V}_{k+2}(K)=\left\{v \in H^2(K):  {\Delta^2 v \in \mathbb{P}_{k-2}(K)}, v|_e \in \mathbb{P}_{k+2}(e),  \partial_{\bm n} v|_e \in \mathbb{P}_{k+1}(e), ~~e \subset \partial K\right\}.
\]
To present the degrees of freedom (DoFs), we introduce a scaled monomial $\mathbb{M}_r(D)$ on a $d$-dimensional domain $D$
\[
\mathbb  M_{r} (D):= \Big \{ \Big ( \frac{\boldsymbol x -  \boldsymbol x_D}{h_D}\Big )^{\boldsymbol  s}, \quad |\boldsymbol  s|\le r\Big \},
\]
where $h_D$ is the diameter of $D$, $\boldsymbol  x_D$ the centroid of $D$, and $r$ a non-negative integer. For the multi-index ${\boldsymbol{s}} \in {\mathbb{N}^d}$, we follow the usual notation
\[\boldsymbol{x}^{\boldsymbol{s}} = x_1^{s_1} \cdots x_d^{s_d},\quad |\boldsymbol{s}| = s_1 +  \cdots  + s_d.\]
Conventionally, $\mathbb  M_r (D) =\{0\}$ for $r\le -1$.
This modified local space can be equipped with the following DoFs (cf. \cite{Brezzi-Marini-2013,Chinosi-Marini-2016}):
\begin{itemize}
	\item  $\tilde{\chi}^v:$ the values of $v$ at the vertices of $K$,
   \[\tilde{\chi}_a^v(v) = v(a), \quad \mbox{$a$ is a vertex of $K$}.\]
	\item  $\tilde{\chi}^g$ : the values of $h_a \nabla v$ at the vertices  of $K$,
  \[\tilde{\chi}_a^g(v) = h_a \nabla v(a), \quad \mbox{$a$ is a vertex of $K$},\]
  where $h_a$ is a characteristic length attached to each vertex $a$, for instance, the average of the diameters of the elements having $a$ as a vertex.
	\item  $\tilde{\chi}^e$ : the moments of $v$ on edges up to degree $k-2$,
  \[ \tilde{\chi}_e(v) = |e|^{-1}(m_e, v)_e, \quad m_e \in \mathbb{M}_{k-2}(e), \quad e \subset\partial K.\]
	\item  $\tilde{\chi}^n$ : the moments of $\partial_{\bm n_e} v$ on edges up to degree $k-1$,
  \[\tilde{\chi}_e^n = (m_e,\partial_{\bm n_e} v)_e, \quad  m_e \in \mathbb{M}_{k-1}(e), \quad e \subset \partial K.\]
	\item  $\tilde{\chi}^K$ : the moments on element $K$ up to degree $k-2$,
  \[\tilde{\chi}_K(v) = |K|^{-1}(m_K,v)_K, \quad m_K \in \mathbb{M}_{k-2}(K).\]
\end{itemize}

Given $v_h \in \widetilde{V}_{k+2}(K)$, the usual definition of the $H^1$-elliptic projection $\Pi_K^{\nabla} v_h \in$ $\mathbb{P}_k(K)$ is described by the following equations:
\begin{equation}\label{H1def}
\left\{\begin{aligned}
 (\nabla \Pi_K^{\nabla} v_h, \nabla q)_K  & = (\nabla v_h, \nabla q)_K, \quad q \in \mathbb{P}_k(K), \\
 \sum\limits_{a \in \mathcal{V}_K} \Pi_K^{\nabla} v_h(a)  & = \sum\limits_{a \in \mathcal{V}_K} v_h(a),
\end{aligned}\right.
\end{equation}
where $\mathcal{V}_K$ is the set of the vertices of $K$. This elliptic projection can be computed by the previous DoFs of $v_h$ by checking the right-hand side of the integration by parts formula:
\[(\nabla v_h, \nabla q)_K = -(v_h, \Delta q)_K + \sum\limits_{e\subset \partial K} \int_e v_h \partial_{\bm n_e} q \mathrm{d} s.\]
However, the goal of the IPVEM is to make $\Pi_K^\nabla v_h$ computable by only using the DoFs of $H^1$-conforming virtual element spaces given by (cf. \cite{Beirao-Brezzi-Cangiani-2013,Ahmad-Alsaedi-Brezzi-2013})
\[V_h^{1,c}(K): = \{ v \in H^1(K): \Delta v|_K \in \mathbb{P}_{k - 2}(K)~~{\text{in}}~~K,\quad  v|_{\partial K} \in \mathbb{B}_k(\partial K) \},\]
where
\[\mathbb{B}_k(\partial K): = \{ v \in C(\partial K):  v|_e \in \mathbb{P}_k(e), \quad e \subset \partial K \}.\]
To do so, Ref.~\cite{ZMZW2023IPVEM} considered the approximation of the RHS by some numerical formula. In view of the $\mathbb{P}_k$ accuracy, namely $\Pi_K^\nabla v_h = v_h$ for $v_h\in \mathbb{P}_k(K)$, the modified $H^1$-projection is defined as
\[
\left\{\begin{aligned}
(\nabla \Pi_K^{\nabla} v_h, \nabla q)_K & = -(v_h, \Delta q)_K + \sum\limits_{e\subset \partial K} Q_{2k-1}^e( v_h \partial_{\bm n_e} q), \quad q \in \mathbb{P}_k(K), \\
	\sum_{a \in \mathcal{V}_K} \Pi_K^{\nabla} v_h(a) & =\sum_{a \in \mathcal{V}_K} v_h(a),
\end{aligned}\right.
\]
with
\[ Q_{2k-1}^e v := |e| \sum\limits_{i=0}^k \omega_i v(\bm{x}_i^e) \approx \int_e v(s) \d s,\]
where $(\omega_i, \bm{x}_i^e)$ are the $(k+1)$ Gauss-Lobatto quadrature weights and points with $\bm{x}_0^e$ and $\bm{x}_k^e$ being the endpoints of $e$.  The piecewise $H^1$-projector $\Pi_h^{\nabla}$ is defined by setting $\Pi_h^{\nabla}|_K=\Pi_K^{\nabla}$ for all $K\in \mathcal{T}_h$.
In this case, the first and the third types of DoFs of $v_h$ on $e$ should be replaced by the $(k+1)$ Gauss-Lobatto points. Notice that when $v_h\in \mathbb{P}_k(K)$, the integrand of $\int_e v_h \partial_{\bm n_e} q \mathrm{d} s$ is a polynomial of order $2k-1$ on $e$, while the algebraic accuracy is exactly $2k-1$ for the $(k+1)$ Gauss-Lobatto points. In particular, we have $2k-1 = 3$ for $k=2$, corresponding to the Simpson's rule.

Following the similar arguments in \cite[Lemma 3.4]{ZMZW2023IPVEM} and \cite[Corollary 3.7]{ZMZW2023IPVEM}, we may derive the inverse inequalities and the boundedness of the $H^1$ projector described in the following Lemma.
\begin{lemma} \label{lem:Pih}
For all $\,v_h\in \widetilde{V}_{k+2}(K)$ there hold
	\begin{align}
		&|v_h|_{m,K}\lesssim h_K^{s-m}|v_h|_{s,K} ,\label{inv}\\
		&|\Pi_K^{\nabla}v_h|_{m,K}\lesssim |v_h|_{m,K},\label{Pi}
	\end{align}
	where $m=1,2$ and $0\le s\le m$.
\end{lemma}

\subsection{The IP virtual element spaces} \label{subsec:IPVEspace}

The $H^1$-elliptic projection is uniquely determined by the DoFs for the $H^1$-conforming virtual element space $V_h^{1,c}$. Because of this, one can replace the additional DoFs of $v \in \widetilde{V}_{k+2}(K)$ by the ones of $\Pi_K^{\nabla} v$. The local interior penalty space is then defined as
\[
V_k(K)= \Big\{v \in \widetilde{V}_{k+2}(K): \tilde{\chi}^g(v)=\tilde{\chi}^g (\Pi_K^{\nabla} v ), \tilde{\chi}^n(v)=\tilde{\chi}^n (\Pi_K^{\nabla} v )\Big\},
\]
which satisfies $\Pi_K^{\nabla} v=v$ for $v \in \mathbb{P}_k(K)$, $\mathbb{P}_k(K) \subset V_k(K)$, and $V_k(K) \subset H^2(K)$.
The associated DoFs are then given by
\begin{itemize}
	\item  $\chi^v:$ the values of $v(a)$, $a \in \mathcal{V}_K$;
	\item $\chi^e:$ the values of $v(\boldsymbol{x}_i^e)$, $i=1,2 \cdots, k-1$, $e \subset \partial K$;
	\item $\chi^K$ : the moments $|K|^{-1}( m_K, v)_K$, $m_K \in \mathbb{M}_{k-2}(K)$.
\end{itemize}
Furthermore, we use $V_h$ to denote the global space of nonconforming virtual element, which is defined piecewise and required that the degrees of freedom are single-valued and vanish for the boundary DoFs.

Since $V_h$ and $V_h^{1,c}$ share the same DoFs, we can introduce the interpolation from $V_h$ to $V_h^{1,c}$.
For any given $v_h\in V_h$, let $I_h^c v_h$ be the nodal interpolant of $v_h$ in $V_h^{1,c}$.
One can define the $H^1$-elliptic projection $\Pi_h^\nabla I_h^c v_h$ as in \eqref{H1def} and find that
\begin{equation}\label{conforminginterp}
	\Pi_h^\nabla I_h^c v_h =  \Pi_h^\nabla v_h, \quad v_h\in V_h
\end{equation}
since $\bm {\chi}(I_h^c v_h) = \bm {\chi}(v_h)$ and $\Pi_h^\nabla$ is uniquely determined by the DoFs in $\bm {\chi}(v_h)$. Here and below, $\bm {\chi}$ is the collection of the DoFs for $V_h$.


As usual, we can define the $H^2$-projection operator $\Pi_K^{\Delta}: V_k(K) \to \mathbb{P}_k(K)$ by finding the solution $\Pi_K^{\Delta} v \in \mathbb{P}_k(K)$ of
\[\begin{cases}
	a^K (\Pi_K^{\Delta} v, q )=a^K(v, q),\quad  q \in \mathbb{P}_k(K), \\
	\widehat{\Pi_K^{\Delta} v}=\widehat{v}, \quad \widehat{\nabla \Pi_K^{\Delta} v}=\widehat{\nabla v}
\end{cases}\]
for any given $v \in V_k(K)$, where the quasi-average $\widehat{v}$ is defined by
\[\widehat{v}=\frac{1}{|\partial K|} \int_{\partial K} v \d s .\]

\subsection{The IP virtual element method with modified jump and penalty terms}
Given the discrete bilinear form
\[
a_h^K(v, w)=a^K(\Pi_K^{\Delta} v, \Pi_K^{\Delta} w)+S^K(v-\Pi_K^{\Delta} v, w-\Pi_K^{\Delta} w), \quad v, w \in V_k(K),
\]
with
\[
S^K(v-\Pi_K^{\Delta} v,w-\Pi_K^{\Delta} w)=h_K^{-2}\bm\chi(v-\Pi_K^{\Delta} v) \cdot  \bm \chi(w-\Pi_K^{\Delta} w).
\]
Let
\[
J_1(v,w)= \sum_{e\in\mathcal{E}_h}\frac{\lambda_e}{|e|}\int_e\Big[\frac{\partial \Pi_h^\nabla v}{\partial \bm{n}_e}\Big]\Big[\frac{\partial \Pi_h^\nabla w}{\partial \bm{n}_e}\Big] \d s,
\]
where  $\lambda_e\ge 1$ is some edge-dependent parameter, and the additional terms are given by
\[
J_2(v,w)=-\sum_{e \in \mathcal{E}_h}\int_e\Big\{\frac{\partial^2\Pi_h^{\nabla}v}{\partial \bm{n}_e^2}\Big\}\Big[\frac{\partial \Pi_h^{\nabla} w}{\partial \bm{n}_e}\Big] \d s,\]
\[J_3(v,w)=-\sum_{e\in\mathcal{E}_h}\int_e\Big\{\frac{\partial^2 \Pi_h^{\nabla}w}{\partial \bm{n}_e^2}\Big\}\Big[\frac{\partial \Pi_h^{\nabla} v}{\partial \bm{n}_e}\Big] \d s. \]
 In contrast to the jump and average terms in \cite{ZMZW2023IPVEM}, here we include the elliptic projector $\Pi_h^\nabla$ for all $v$ and $w$ in $J_i$ for $i=1,2,3$ to simplify the implementation. According to the definition of $V_k(K)$, one easily finds that $J_2$ (or $J_3$) defined here coincides with the one given in \cite{ZMZW2023IPVEM}, viz.,
\[J_2(v,w)=-\sum_{e \in \mathcal{E}_h}\int_e\Big\{\frac{\partial^2\Pi_h^{\nabla}v}{\partial \bm n_e^2}\Big\}\Big[\frac{\partial \Pi_h^{\nabla}w}{\partial \bm n_e}\Big] \d s = -\sum_{e \in \mathcal{E}_h}\int_e\Big\{\frac{\partial^2\Pi_h^{\nabla}v}{\partial \bm n_e^2}\Big\}\Big[\frac{\partial w}{\partial \bm n_e}\Big] \d s.\]

The IPVEM for the problem \eqref{originalSP} can be described as follows: Find $u_h\in  V_h$ such that
\begin{equation}\label{IPVEM1}
	\varepsilon^2\mathcal{A}_h(u_h,v_h)+b_h(u_h, v_h) = \langle f_h,v_h\rangle,\quad v_h\in V_h,
\end{equation}
where
\begin{equation}\label{BilinearAh}
	\mathcal{A}_h(u_h,v_h) = a_h(u_h,v_h)+J_1(u_h,v_h)+J_2(u_h,v_h)+J_3(u_h,v_h),
\end{equation}
and $b_h(v,w)=\sum\limits_{K\in \mathcal{T}_h} b_h^K(v,w)$, with
\[b_h^K(v,w):=(\nabla \Pi_K^{\nabla}v,\nabla \Pi_K^{\nabla}w)+\bm \chi(v-\Pi_K^{\nabla} v)\cdot \bm \chi (w-\Pi_K^{\nabla} w).\]
The right-hand side is
\[
\langle f_h,v_h \rangle:	=
\sum\limits_{K\in \mathcal{T}_h}\int_{K} f\,\Pi_{0,K}^k v_h \d x,
\]	
with $\Pi_{0,K}^k$ being the $L^2$ projector onto $\mathbb{P}_k(K)$.

\begin{remark} \label{rem:observation}
In view of Eq.~\eqref{conforminginterp}, we have
\begin{equation}
	b_h^K(v,w) = b_h^K(v, I_h^c w), \quad v, w \in V_h,\label{bIh}
\end{equation}
which is crucial in the error analysis.
We also have
\begin{equation}
b_h^K(I_h^cv,I_h^cv)=b_h^K(v,I_h^cv), \quad v, w \in V_h,
\end{equation}
which implies
\begin{equation}\label{boundIhc}
|I_h^c v_h|_{1,h}\lesssim |v_h|_{1,h}.
\end{equation}
\end{remark}

In what follows, we define
\begin{equation} \label{normm}
	\|w\|_h^2:=|w|_{2,h}^2+J_1(w,w), \qquad
	\| w \|_{\varepsilon,h} ^2:=\varepsilon^2\|w\|_h^2+|w|_{1,h}^2.
\end{equation}

\begin{lemma}
$\|\cdot\|_h$ and $\| \cdot \|_{\varepsilon,h}$ are norms on $V_h$.
\end{lemma}
\begin{proof}
It is enough to prove that $\|v_h\|_h = 0$ implies $v_h = 0$ for any given $v_h \in V_h$. By definition, $\|v_h\|_h = 0$ is equivalent to
\[|v_h|_{2,h} = 0 , \qquad J_1(v_h,v_h) = 0.\]
The equation $|v_h|_{2,h} = 0$ shows that $\nabla_h v_h$ is a piecewise constant on $\mathcal{T}_h$. On the other hand, the direct manipulation yields
 \begin{align}\label{eq:1}
\int_e\Big[\frac{\partial \Pi_h^\nabla v}{\partial  \bm{n}_e}\Big]\Big[\frac{\partial \Pi_h^\nabla w}{\partial  \bm{n}_e}\Big] \d s
& =: \int_e q_{k-1}^e  \Big[\frac{\partial \Pi_h^\nabla w}{\partial  \bm{n}_e}\Big] \d s
  = \int_e q_{k-1}^e  \Big[\frac{\partial w}{\partial  \bm{n}_e}\Big] \d s\notag\\
& = \int_e q_{k-1}^e  \Pi^{k-1}_{0,e} \Big[\frac{\partial w}{\partial  \bm{n}_e}\Big] \d s
  = \int_e \Big[\frac{\partial \Pi_h^\nabla v}{\partial  \bm{n}_e}\Big]  \Pi^{k-1}_{0,e} \Big[\frac{\partial w}{\partial  \bm{n}_e}\Big] \d s \notag \\
& = \int_e \Big[\frac{\partial v}{\partial  \bm{n}_e}\Big]  \Pi^{k-1}_{0,e} \Big[\frac{\partial w}{\partial  \bm{n}_e}\Big]  \d s
  = \int_e \Pi^{k-1}_{0,e} \Big[\frac{\partial v}{\partial  \bm{n}_e}\Big]  \Pi^{k-1}_{0,e} \Big[\frac{\partial w}{\partial  \bm{n}_e}\Big]  \d s ,
\end{align}
where $v, w \in V_h$.
Thus, $J_1(v_h,v_h) = 0$ implies $\Pi^{k-1}_{0,e} \Big[\frac{\partial v_h}{\partial  \bm{n}_e}\Big] = 0$, where $k \ge 2$. Since $\nabla_h v_h$ is piecewise constant, we further obtain $\Big[\frac{\partial v_h}{\partial  \bm{n}_e}\Big] = 0$ over the edges of $\mathcal{T}_h$. That is, the normal derivative of $v_h$ is continuous over interior edges and vanishes on the boundary of $\Omega$. This reduces to the discussion in the proof of Lemma 4.2 in \cite{ZMZW2023IPVEM}, so we omit the remaining argument.
\end{proof}

\subsection{Well-posedness of the discrete problem}

In \cite{ZMZW2023IPVEM}, the mesh-dependent parameter $\lambda_e$ was chosen as a sufficiently large constant $\lambda$. The stability for the bilinear forms can be obtained by using the similar arguments used in the proof of Theorem 4.3 in \cite{ZMZW2023IPVEM}. We omit the details with the results described as follows.
\begin{itemize}
  \item[-] $k$-consistency: for all $v \in V_k(K)$ and $q \in \mathbb{P}_k(K)$, it holds that
\begin{equation}\label{consistent}
a_h^K(v, q)=a^K(v, q), \quad b_h^K(v, q)=b^K(v, q).
\end{equation}
  \item[-] Stability: there exist two positive constants $\alpha_*$ and $\alpha^*$, independent of $h$, such that
\begin{align}
  & \alpha_*a^K(v, v) \le a_h^K(v, v) \le  \alpha^* a^K(v, v), \label{stable} \\
  & \alpha_* b^K(v,v) \le  b_h^K(v,v) \le \alpha^* b^K(v,v) \label{stableb}
\end{align}
  for all $v \in V_k(K)$.
\end{itemize}

Here, we aim to provide an automated mesh-dependent selection of $\lambda_e$ following a similar deduction with that described in \cite{Carsten2023IP}.

\begin{lemma} \label{lem:parameter}
For every constant $a$ satisfying
\[ a> \max \Big\{ 1,  \frac{1}{\sqrt{\alpha_*}}\Big\},\]
where $\alpha_*$ is given in \eqref{stable}, define the penalty parameter as
\[\lambda_e = \begin{cases}
\dfrac{a N_K k(k-1) h_e^2 }{4} \Big( \dfrac{1}{|T_+|} + \dfrac{1}{|T_-|} \Big),  & \qquad e \in \mathcal{E}^0,\\
\dfrac{a N_K k(k-1) h_e^2 }{2 |T_+|} ,  &  \qquad e \in \mathcal{E}^{\partial},
\end{cases}\]
where $T_+$ and $T_-$ are the neighboring virtual triangles for an interior edge $e$, $T_+$ is the adjacent virtual triangle for a boundary edge $e$, and $N_K$ is the maximum number of edges of elements.
Then there holds
\[\mathcal{A}_h(v_h,v_h) \ge \kappa \|v_h\|_h^2\]
with the constant
\[\kappa = \min \{\alpha_*, 1\} - \frac{1}{\sqrt{a}}.\]
\end{lemma}

\begin{proof}
For every $\kappa>0$, consider the difference
\begin{align*}
\mathcal{A}_h(v_h,v_h) - \kappa \|v_h\|_h^2
& = ( a_h(v_h,v_h) - \kappa |v_h|_{2,h}^2 ) + (1-\kappa) J_1(v_h,v_h)+J_2(v_h,v_h)+J_3(v_h,v_h) \\
& = ( a_h(v_h,v_h) - \kappa |v_h|_{2,h}^2 ) + (1-\kappa) J_1(v_h,v_h) + 2 J_2(v_h,v_h) \\
& \ge (\alpha_* - \kappa) |v_h|_{2,h}^2 + (1-\kappa) J_1(v_h,v_h) + 2 J_2(v_h,v_h),
\end{align*}
where the stability \eqref{stable} is used.
For every $\epsilon>0$, the Young's inequality gives
\begin{align*}
 2 |J_2(v_h,v_h)|
& = 2 \Big | \sum_{e \in \mathcal{E}_h}\int_e\Big\{\frac{\partial^2\Pi_h^{\nabla}v_h}{\partial \bm{n}_e^2}\Big\}\Big[\frac{\partial \Pi_h^{\nabla} v_h}{\partial \bm{n}_e}\Big] \d s \Big| \\
& \le \sum_{e \in \mathcal{E}_h} \Big(
        \frac{\epsilon \lambda_e}{h_e} \Big\| \Big[\frac{\partial \Pi_h^{\nabla} v_h}{\partial \bm{n}_e}\Big] \Big\|_{0,e}^2
      + \frac{h_e}{\epsilon \lambda_e} \Big\| \Big\{\frac{\partial^2\Pi_h^{\nabla}v_h}{\partial \bm{n}_e^2}\Big\} \Big\|_{0,e}^2\Big) \\
& = \epsilon J_1(v_h, v_h) + \sum_{e \in \mathcal{E}_h} \frac{h_e}{\epsilon \lambda_e} \Big\| \Big\{\frac{\partial^2\Pi_h^{\nabla}v_h}{\partial \bm{n}_e^2}\Big\} \Big\|_{0,e}^2.
\end{align*}
This reduces to the estimate of the average $\Big\{\frac{\partial^2\Pi_h^{\nabla}v_h}{\partial \bm{n}_e^2}\Big\} = \Big\{\bm{n}_e^T (\nabla_h^2\Pi_h^{\nabla}v_h) \bm{n}_e\Big\}$.

For an interior edge $e$ with the neighboring virtual triangles $T_+$ and $T_-$, the definition of the average and the Young's inequality give
\begin{align}
\Big\| \Big\{\frac{\partial^2\Pi_h^{\nabla}v_h}{\partial \bm{n}_e^2}\Big\} \Big\|_{0,e}^2
& = \Big\| \frac12 n_e^T ( (\nabla_h^2\Pi_h^{\nabla}v_h)|_{T_+} + (\nabla_h^2\Pi_h^{\nabla}v_h)|_{T_-}) n_e \Big\|_{0,e}^2  \nonumber \\
& \le \frac{1}{4} \| (\nabla_h^2\Pi_h^{\nabla}v_h)|_{T_+} +  (\nabla_h^2\Pi_h^{\nabla}v_h)|_{T_-} \|_{0,e}^2 \nonumber \\
& = \frac{1}{4}(  \| (\nabla_h^2\Pi_h^{\nabla}v_h)|_{T_+} \|_{0,e}^2 +  \|(\nabla_h^2\Pi_h^{\nabla}v_h)|_{T_-} \|_{0,e}^2
   + 2 \int_e (\nabla_h^2\Pi_h^{\nabla}v_h)|_{T_+}(\nabla_h^2\Pi_h^{\nabla}v_h)|_{T_-} \d s \nonumber \\
& \le \frac14 \Big( (1+\alpha) \| (\nabla_h^2\Pi_h^{\nabla}v_h)|_{T_+} \|_{0,e}^2 + (1+\frac{1}{\alpha})\|(\nabla_h^2\Pi_h^{\nabla}v_h)|_{T_-} \|_{0,e}^2  \Big) \label{nnvh}
\end{align}
for any $\alpha>0$. Since $(\nabla_h^2\Pi_h^{\nabla}v_h)|_{T_\pm}$ are polynomials of degree $k-2$,  this and Lemma 3.4 in \cite{Carsten2023IP} result in
\[\Big\| \Big\{\frac{\partial^2\Pi_h^{\nabla}v_h}{\partial \bm{n}_e^2}\Big\} \Big\|_{0,e}^2
\le \frac{k(k-1)h_e}{8}  \Big( \frac{1+\alpha}{|T_+|} \| \nabla_h^2\Pi_h^{\nabla}v_h \|_{0,T_+}^2 + \frac{1+1/\alpha}{|T_-|}\|\nabla_h^2\Pi_h^{\nabla}v_h \|_{0,T_-}^2  \Big) .\]
The optimal value $\alpha = |T_+|/|T_-|$ and the boundedss of $\Pi_h^{\nabla}$ in \eqref{Pi} of Lemma \ref{lem:Pih} lead to
\begin{align}
\Big\| \Big\{\frac{\partial^2\Pi_h^{\nabla}v_h}{\partial \bm{n}_e^2}\Big\} \Big\|_{0,e}^2
& \le \frac{k(k-1)h_e}{8}  \Big( \frac{1}{|T_+|}  + \frac{1}{|T_-|} \Big) \| \nabla_h^2\Pi_h^{\nabla}v_h \|_{0,T_+ \cup T_-} \nonumber\\
& \le \frac{k(k-1)h_e}{8}  \Big( \frac{1}{|T_+|}  + \frac{1}{|T_-|} \Big) | v_h |_{2,K_+ \cup K_-}. \label{average}
\end{align}
For a boundary edge $e$ with adjacent virtual triangle $T_+$, we can get
\[\Big\| \Big\{\frac{\partial^2\Pi_h^{\nabla}v_h}{\partial \bm{n}_e^2}\Big\} \Big\|_{0,e}^2
\le \frac{k(k-1)h_e}{2 |T_+|}  \| \nabla_h^2\Pi_h^{\nabla}v_h \|_{0,T_+}
\le \frac{k(k-1)h_e}{2 |T_+|}  |v_h |_{2,K_+} .\]
Consequently, with the choice of $\lambda_e$, the sum over all edges and the finite overlap of the edge patches lead to
\begin{equation*}
2 |J_2(v_h,v_h)| \le \epsilon J_1(v_h, v_h) + \frac{1}{a\epsilon} |v_h|_{2,h}.
\end{equation*}

The above discussion gives
\begin{align*}
\mathcal{A}_h(v_h,v_h) - \kappa \|v_h\|_h^2
& \ge (\alpha_* - \frac{1}{a\epsilon}- \kappa) |v_h|_{2,h}^2 + (1 - \epsilon-\kappa) J_1(v_h,v_h).
\end{align*}
Every choice of $0< \kappa \le \min\{ \alpha_* - 1/(a\epsilon) ,  1 - \epsilon\}$ leads to a nonnegative lower bound, and so proves the claim that
\[\mathcal{A}_h(v_h,v_h) \ge \kappa \|v_h\|_h^2.\]
Taking $\epsilon = 1/\sqrt{a}$ with $ a> \max \{ 1,  \frac{1}{\sqrt{\alpha_*}}\}$ results in
\[\kappa := \min \Big\{ \alpha_* - \frac{1}{\sqrt{a}} ,  1 - \frac{1}{\sqrt{a}} \Big\}
= \min \{\alpha_*, 1\} - \frac{1}{\sqrt{a}},\]
as required.
\end{proof}

\begin{remark}
We remark that the lower and upper bounds of $\lambda_e$ are independent of mesh sizes, namely, $\lambda_e \eqsim 1$, under the given mesh assumption.
\end{remark}

\begin{theorem} For the parameter $\lambda_e$ chosen as in Lemma \ref{lem:parameter}, there exists a unique solution to the discrete problem \eqref{IPVEM1}.
\end{theorem}
\begin{proof}
For any $v_h \in V_h$, Lemma \ref{lem:parameter} along with the stability \eqref{stableb} yields the coercivity
\begin{equation}\label{coercivityEps}
\varepsilon^2\mathcal{A}_h(v_h,v_h)+b_h(v_h,v_h) \ge \varepsilon^2 \kappa \|v_h\|_h^2 + \alpha_* |v_h|_{1,h}^2
\ge \min \{\kappa, \alpha_*\} \| v_h \|_{\varepsilon,h}^2 = \alpha_*\| v_h \|_{\varepsilon,h}^2
\end{equation}
For the boundedness of the bilinear form, by the definition of $\normm{\cdot}_h$ given in \eqref{normm}, it suffices to consider $\varepsilon^2 J_2(v_h,w_h)$ for any $v_h,w_h \in V_h$. The Cauchy-Schwarz inequality gives
\begin{align}
|J_2(v_h,w_h)|
& = \Big|\sum_{e \in \mathcal{E}_h} \int_e\Big\{\frac{\partial^2\Pi_h^{\nabla}v_h}{\partial \bm{n}_e^2}\Big\}\Big[\frac{\partial \Pi_h^{\nabla} w_h}{\partial \bm{n}_e}\Big] \d s\Big| \nonumber \\
& \lesssim \Big(\sum_{e \in \mathcal{E}_h}\frac{h_e}{\lambda_e} \Big\|\Big\{\frac{\partial^2\Pi_h^{\nabla}v_h}{\partial \bm{n}_e^2}\Big\}\Big\|_e^2\Big)^{\frac{1}{2}}\Big(\sum_{e \in \mathcal{E}_h}\frac{\lambda_e}{h_e}\Big\|\Big[\frac{\partial \Pi_h^{\nabla} w_h}{\partial \bm{n}_e}\Big]\Big\|_e^2\Big)^{\frac{1}{2}}  \nonumber \\
& = \Big(\sum_{e \in \mathcal{E}_h}\frac{h_e}{\lambda_e} \Big\|\Big\{\frac{\partial^2\Pi_h^{\nabla}v_h}{\partial \bm{n}_e^2}\Big\}\Big\|_e^2\Big)^{\frac{1}{2}} J_1(w_h,w_h)^{1/2}. \label{J2bound}
\end{align}
For an interior edge $e$ with the neighboring virtual triangles $T_+$ and $T_-$, we obtain from \eqref{average} and the definition of $\lambda_e$ that
\begin{align*}
\frac{h_e}{\lambda_e} \Big\| \Big\{\frac{\partial^2\Pi_h^{\nabla}v_h}{\partial \bm{n}_e^2}\Big\} \Big\|_{0,e}^2
 \le \frac{k(k-1)h_e}{8}  \Big( \frac{1}{|T_+|}  + \frac{1}{|T_-|} \Big) | v_h |_{2,K_+ \cup K_-}
 = \frac{1}{2aN_K}   | v_h |_{2,K_+ \cup K_-},
\end{align*}
where $a$ and $N_K$ are defined in Lemma \ref{lem:parameter}. For a boundary edge $e$, one can get
\begin{align*}
\frac{h_e}{\lambda_e} \Big\| \Big\{\frac{\partial^2\Pi_h^{\nabla}v_h}{\partial \bm{n}_e^2}\Big\} \Big\|_{0,e}^2
 \le  \frac{1}{aN_K}   | v_h |_{2,K_+}.
\end{align*}
Therefore,
\[\Big(\sum_{e \in \mathcal{E}_h}\frac{h_e}{\lambda_e} \Big\|\Big\{\frac{\partial^2\Pi_h^{\nabla}v_h}{\partial \bm{n}_e^2}\Big\}\Big\|_e^2\Big)^{\frac{1}{2}} \lesssim |v_h|_{2,h}, \]
which together with \eqref{J2bound} yields
\[|\varepsilon^2 J_2(v_h,w_h)| \lesssim \| v_h \|_{\varepsilon,h} \| w_h \|_{\varepsilon,h}.\]
The proof is finished by using the Lax-Milgram lemma.
\end{proof}

\section{Error analysis} \label{sec:err}

\subsection{An abstract Strang-type lemma}

For any function $v \in H_0^1(\Omega)\cap H^2(\Omega)$, its interpolation $I_h v$ is defined by the condition that $v$ and $I_h v$ have the same degrees of freedom:
\[\chi_i(I_h v) = \chi_i(v), \qquad \chi_i \in \bm {\chi}.\]
Since the computation of the elliptic projection $\Pi_h^\nabla I_h u$ only involves the DoFs of $u$ and $\bm {\chi}(I_h u) = \bm {\chi}(u)$,  we define $\Pi_h^\nabla u= \Pi_h^\nabla I_h u$ in what follows,  which makes the expressions
$\|u-u_h\|_{\varepsilon,h}$, $\|u-I_hu\|_{\varepsilon,h}$ and $\|u-u_{\pi}\|_{\varepsilon,h}$ well-defined in the following lemma, where $u$ is the exact solution to \eqref{originalSP}.

\begin{lemma}\label{lem:StrangIPVEM}
	Assume that $u\in H_0^2(\Omega)\cap H^3(\Omega)$ is the exact solution to \eqref{originalSP}. Then we have the error estimate
\begin{equation}\label{Strangerr}
\|u-u_h\|_{\varepsilon,h}\lesssim \|u-I_hu\|_{\varepsilon,h} + \| u-u_{\pi} \|_{\varepsilon,h}  + \|f-f_h\|_{V_h'} + E_h,
\end{equation}
	where $u_\pi$ is a piecewise polynomial,
	\[\|f-f_h\|_{V_h'} = \sup _{v_h \in V_h\backslash \{0\}} \frac{\langle f - f_h, v_h\rangle}{\| v_h \|_{\varepsilon,h}}\]
	and the consistency term
	\begin{equation}\label{consistencyterm}
		E_h = \sup _{v_h \in V_h \backslash\{0\}} \frac{E_A(u,v_h) + E_J(u,v_h)}{\|v_h\|_{\varepsilon,h}},
	\end{equation}
	with
	\[E_A(u,v_h) = \sum\limits_{K\in \mathcal{T}_h}(\varepsilon^2 a^K(u, v_h)+b^K(u,I_h^cv_h)) - (f, v_h),\]
	\[E_J(u,v_h) = \varepsilon^2(J_1(I_h u, v_h) + J_2(I_h u, v_h) + J_3(I_h u, v_h)),\]
where $I_h^c$ is the conforming interpolation operator defined in Section \ref{subsec:IPVEspace}.
\end{lemma}
\begin{proof}
Let $\delta_h:= I_hu-u_h$. We have by the coercivity \eqref{coercivityEps} and the definition of the discrete problem \eqref{IPVEM1} that
\begin{align*}
\alpha_*\|\delta_h\|_{\varepsilon,h}^2
& \le \varepsilon^2 \Big( a_h(\delta_h, \delta_h)+J_1(\delta_h, \delta_h)+J_2(\delta_h, \delta_h) +J_3(\delta_h, \delta_h)\Big)+b_h(\delta_h,\delta_h) \\
& = \varepsilon^2 \Big( a_h(I_h u, \delta_h)+J_1(I_h u, \delta_h)+J_2(I_h u, \delta_h)+J_3(I_h u, \delta_h)\Big)+ b_h(I_hu,\delta_h)-\langle f_h, \delta_h\rangle \\
& = \varepsilon^2 a_h(I_h u, \delta_h)+ b_h(I_hu,I_h^c\delta_h)
		- (f, \delta_h) + \langle f - f_h, \delta_h\rangle + E_J(u,v_h),
\end{align*}
where we have used \eqref{bIh} in Remark \ref{rem:observation}.

According to  the $k$-consistency \eqref{consistent} and the stability formulas \eqref{stable} and \eqref{stableb}, combing with \eqref{boundIhc}, we get
	\begin{align*}
		& \varepsilon^2 a_h(I_h u, \delta_h)+b_h(I_h u,\delta_h) - (f, \delta_h) \\
		& = \sum_{K \in \mathcal{T}_h}\varepsilon^2 (a_h^K(I_h u-u_\pi, \delta_h)+a^K(u_\pi-u, \delta_h))
		+\sum_{K\in \mathcal{T}_h}(b_h^K(I_h u-u_\pi, \delta_h)+b^K(u_\pi-u, I_h^c\delta_h)) \\
		& \qquad + \sum\limits_{K\in\mathcal{T}_h}( \varepsilon^2 a^K(u,\delta_h) + b^K(u,I_h^c\delta_h) ) - (f, \delta_h)  \\
		& \le \varepsilon^2 (|u-I_h u|_{2,h} + |u-u_\pi|_{2,h})|\delta_h|_{2,h}+(|u-I_h u|_{1,h} + |u-u_\pi|_{1,h})|\delta_h|_{1,h}
		+ E_A(u,v_h)  \\
		& \lesssim (\| u-I_h u \|_{\varepsilon,h} + \| u-u_\pi \|_{\varepsilon,h}) \| \delta_h \|_{\varepsilon,h}
		+ E_A(u,v_h).
	\end{align*}
	The proof is completed by using the triangle inequality.
\end{proof}

\begin{lemma}\label{rem:Strangdecomp}
One can get rid of $\varepsilon^2 J_2(I_hu,v_h)$ in the consistency term \eqref{consistencyterm} and obtain
\begin{align*}
E_A(u,v_h)+E_J(u,v_h)
 = E_{A1} + E_{A2}+E_{A3}
  + \varepsilon^2(J_1(I_hu,v_h)+J_3(I_hu,v_h)),
\end{align*}
where
\[E_{A1} = \varepsilon^2\sum_{e \in \mathcal{E}_h}\int_e  \Big\{\frac{\partial^2 (u-\Pi_h ^\nabla I_h u)}{\partial \boldsymbol{n}_e^2}\Big\} \left[\frac{\partial  v_h  }{\partial \boldsymbol{n}_e}\right]  \d s, \qquad
E_{A2} = \varepsilon^2\sum_{e \in \mathcal{E}_h}\int_e \frac{\partial^2 u}{\partial \boldsymbol{n}_e \partial \boldsymbol{t}_e}\left[\frac{\partial v_h}{\partial \boldsymbol{t}_e}\right] \d s, \]
	\[E_{A3} = -\varepsilon^2 \sum_{K\in \mathcal{T}_h}(\nabla \Delta u,\nabla v_h)_K +  \sum_{K\in \mathcal{T}_h}(\nabla u,  \nabla I_h^c v_h)_K - (f, v_h).\]
\end{lemma}
\begin{proof}
For brevity, we use the summation convention whereby summation is implied when an index is repeated exactly two times.
For any given $v_h \in V_h$, the integration by parts gives
	\begin{align*}
		a^K(u,v_h)
		= (\nabla^2u, \nabla^2v_h)_K = \int_K \partial_{ij} u \partial_{ij} v_h \d x
		= \int_{\partial K}  \partial_{ij} u n_j\partial_i v_h \d s - \int_K \partial_{ijj} u \partial_i v_h \d x,
	\end{align*}
	where $\bm {n}_K = (n_1, n_2)$. Since $\partial_iv = n_i\partial_{\boldsymbol n}v + t_i\partial_{\boldsymbol t}v$, we immediately obtain
	\begin{align*}
		a^K(u,v_h)
		& = \int_{\partial K}  \partial_{ij} u n_i n_j\partial_{\boldsymbol n}v_h  \d s
		+ \int_{\partial K}  \partial_{ij} u t_i n_j \partial_{\boldsymbol t}v_h \d s - \int_K \partial_{ijj} u \partial_i v_h \d x \\
		& = \int_{\partial K} \frac{\partial^2 u}{\partial \boldsymbol{n}_K^2} \frac{\partial v_h}{\partial \boldsymbol{n}_K} \d s+\int_{\partial K} \frac{\partial^2 u}{\partial \boldsymbol{n}_K \partial \boldsymbol{t}_K} \frac{\partial v_h}{\partial \boldsymbol{t}_K} \d s- (\nabla \Delta u,\nabla v_h)_K.
	\end{align*}
Therefore,
	\begin{align}
		E_A(u,v_h)
		& : = \sum\limits_{K\in \mathcal{T}_h}(\varepsilon^2 a^K(u, v_h)+b^K(u,I_h^cv_h)) - (f, v_h) \nonumber\\
		& = \varepsilon^2\sum_{e \in \mathcal{E}_h}\int_e \Big(\frac{\partial^2 u}{\partial \boldsymbol{n}_e^2}\left[\frac{\partial v_h}{\partial \boldsymbol{n}_e}\right] + \frac{\partial^2 u}{\partial \boldsymbol{n}_e \partial \boldsymbol{t}_e}\left[\frac{\partial v_h}{\partial \boldsymbol{t}_e}\right]\Big) \mathrm{d} s \nonumber \\
		& \qquad -\varepsilon^2 \sum_{K\in \mathcal{T}_h}(\nabla \Delta u,\nabla v_h)_K +  \sum_{K\in \mathcal{T}_h}(\nabla u,  \nabla I_h^cv_h)_K - (f, v_h) \nonumber \\
		& =:  \widetilde{E}_{A1} + E_{A2} + E_{A3}, \label{EA}
	\end{align}
	where
	\[\widetilde{E}_{A1} = \varepsilon^2\sum_{e \in \mathcal{E}_h}\int_e  \frac{\partial^2 u}{\partial \boldsymbol{n}_e^2}\left[\frac{\partial v_h}{\partial \boldsymbol{n}_e}\right] \d s,\]
   and $E_{A2}$ and $E_{A3}$ are defined as in the lemma.
For $\widetilde{E}_{A1}$, we consider the decomposition
	\begin{align*}
	\widetilde{E}_{A1}
		& = \varepsilon^2\sum_{e \in \mathcal{E}_h}\int_e \Big\{ \frac{\partial^2 (u-\Pi_h ^\nabla I_h u)}{\partial \boldsymbol{n}_e^2}\Big\}\left[\frac{\partial v_h}{\partial \boldsymbol{n}_e}\right] \d s +\varepsilon^2\sum_{e \in \mathcal{E}_h}\int_e  \Big\{\frac{\partial^2 \Pi_h ^\nabla I_h u}{\partial \boldsymbol{n}_e^2}\Big\}\left[\frac{\partial v_h}{\partial \boldsymbol{n}_e}\right] \d s \\
		& = \varepsilon^2\sum_{e \in \mathcal{E}_h}\int_e \Big\{ \frac{\partial^2 (u-\Pi_h ^\nabla I_h u)}{\partial \boldsymbol{n}_e^2}\Big\}\left[\frac{\partial v_h}{\partial \boldsymbol{n}_e}\right] \d s - \varepsilon^2 J_2 (I_h u, v_h)
   = E_{A1} - \varepsilon^2 J_2 (I_h u, v_h).
	\end{align*}
The result follows from Lemma \ref{lem:StrangIPVEM}.
\end{proof}

\subsection{Error estimate}

For $3\le s \le k+1$, the estimate of the load term is
\begin{align*}
\langle f-f_h, v_h \rangle
& = \sum_{K \in \mathcal{T}_h} (f, v_h-\Pi_{0,K}^k v_h )_K =\sum_{K \in \mathcal{T}_h} (f-\Pi_{0,K}^k f, v_h-\Pi_{0,K}^k v_h )_K \nonumber\\
&  \lesssim h^{s-3} \|f \|_{s-3} h | v_h |_{1,h} \lesssim h^{s-2} \|f \|_{s-3} \| v_h \|_{\varepsilon,h},
\end{align*}
namely,
\begin{equation} \label{errRhs}
\| f - f_h\|_{V_h'} \lesssim h^{s-2} \|f \|_{s-3}.
\end{equation}

As shown in Lemma 3.11 of \cite{ZMZW2023IPVEM}, we can derive the interpolation error estimate for our IPVEM.
\begin{lemma}\label{lem:interpIPVEM}
	For any $v\in H_0^1(\Omega)\cap H^s(\Omega)$ with $2\le s\le k+1$, it holds
	\begin{equation*}
		|v-I_hv|_{m,K}\lesssim h_K^{s-m}|v|_{s,K},\quad m=0,1,2,\quad  K\in \mathcal{T}_h.
	\end{equation*}	
\end{lemma}

We now consider the consistency term in Lemma \ref{lem:StrangIPVEM}.
\begin{lemma}\label{errConsistency}
 Assume that $u \in H_0^2(\Omega) \cap H^s(\Omega)$ is the exact solution to \eqref{originalSP} with $3 \leq s \leq k+1$. Then the consistency error is bounded by
	\begin{equation}\label{aim1}
			E_h \lesssim h^{s-2} (\varepsilon |u|_s + \|f\|_{s-3}) .
	\end{equation}

\end{lemma}
\begin{proof} For clarity, we divide the proof into two steps.

Step 1: According to Lemma \ref{rem:Strangdecomp}, one has
\begin{align}\label{eq:mid}
 E_A(u,v_h)+E_J(u,v_h)
& =  E_{A1} + E_{A2}+E_{A3} +\varepsilon^2(J_1(I_hu,v_h)+J_3(I_hu,v_h))\nonumber\\
& = \varepsilon^2\sum_{e \in \mathcal{E}_h}\int_e  \Big\{\frac{\partial^2 (u-\Pi_h ^\nabla I_h u)}{\partial \boldsymbol{n}_e^2}\Big\} \left[\frac{\partial (v_h - \Pi_h^\nabla v_h)}{\partial \boldsymbol{n}_e}\right]  \d s \nonumber \\
& \quad + \varepsilon^2\sum_{e \in \mathcal{E}_h}\int_e  \Big\{\frac{\partial^2 (u-\Pi_h ^\nabla I_h u)}{\partial \boldsymbol{n}_e^2}\Big\} \left[\frac{\partial \Pi_h^\nabla v_h}{\partial \boldsymbol{n}_e}\right]  \d s\nonumber\\
& \quad + E_{A2}+E_{A3}+\varepsilon^2(J_1(I_hu,v_h)+J_3(I_hu,v_h)) \nonumber\\
& \quad =: (E_{A1}^1 + E_{A1}^2) + E_{A2}+E_{A3}+\varepsilon^2(J_1(I_hu,v_h)+J_3(I_hu,v_h)).
\end{align}

Applying the Cauchy-Schwarz inequality, the trace inequality and the projection error estimates, we can get
	\begin{align*}
	E_{A1}^1	
   & = \sum_{e \in \mathcal{E}_h}\int_e  \Big\{\frac{\partial^2 (u-\Pi_h ^\nabla I_h u)}{\partial \boldsymbol{n}_e^2}\Big\} \left[\frac{\partial (v_h - \Pi_h^\nabla v_h)}{\partial \boldsymbol{n}_e}\right]  \d s \\
		& \le \Big\|\Big\{\frac{\partial^2 (u-\Pi_h ^\nabla I_h u)}{\partial \boldsymbol{n}_e^2}\Big\}\Big\|_{0,e} \Big\|\Big[\frac{\partial (v_h - \Pi_h^\nabla v_h)}{\partial \boldsymbol{n}_e}\Big]\Big\|_{0,e} \\
		& \lesssim \sum\limits_{K\in \mathcal{T}_h} \Big(h_K^{1/2} |u-\Pi_h^\nabla I_h u|_{3,K} + h_K^{-1/2}|u-\Pi_h^\nabla I_h u|_{2,K}\Big) \\
		& \qquad \times  \Big(h_K^{1/2} |v_h - \Pi_h^\nabla v_h |_{2,K} + h_K^{-1/2}|v_h - \Pi_h^\nabla v_h |_{1,K}\Big) \\
		& \lesssim \sum\limits_{K\in \mathcal{T}_h} \Big(h_K|u-\Pi_h^\nabla I_h u|_{3,K} + |u-\Pi_h^\nabla I_h u|_{2,K}\Big) |v_h |_{2,K} .
	\end{align*}
	The inverse inequality for polynomials and the boundedness of $\Pi_h^\nabla$  in \eqref{Pi} imply
	\begin{align*}
		& h_K|u-\Pi_h^\nabla I_h u|_{3,K} + |u-\Pi_h^\nabla I_h u|_{2,K} \\
		& \le h_K|u-u_\pi|_{3,K} + h_K|\Pi_h^\nabla (I_h u-u_\pi)|_{3,K}+ |u-u_\pi|_{2,K} + |\Pi_h^\nabla (I_h u-u_\pi)|_{2,K} \\
		& \lesssim h_K|u-u_\pi|_{3,K} + |I_h u-u_\pi|_{2,K}+ |u-u_\pi|_{2,K} + |I_h u-u_\pi|_{2,K} \\
		& \le |u -I_h u|_{2,K} + |u-u_\pi|_{2,K} + h_K|u-u_\pi|_{3,K}.
	\end{align*}
This together with the Dupont-Scott theory and Lemma \ref{lem:interpIPVEM} gives
\begin{align*}
E_{A1}^1	
& \lesssim \varepsilon  (|u -I_h u|_{2,h} + |u-u_\pi|_{2,h} + h |u-u_\pi|_{3,h}) \| v_h \|_{\varepsilon,h}\\
& \lesssim \varepsilon h^{s-2}|u|_s\| v_h \|_{\varepsilon,h}.
\end{align*}
According to the Cauchy-Schwarz inequality,
\begin{align*}
E_{A1}^2
& =	\sum_{e \in \mathcal{E}_h}\int_e  \Big\{\frac{\partial^2 (u-\Pi_h ^\nabla I_h u)}{\partial \boldsymbol{n}_e^2}\Big\} \Big[\frac{\partial \Pi_h^\nabla v_h}{\partial \boldsymbol{n}_e}\Big]  \d s \\
& \le \Big(\sum_{e \in \mathcal{E}_h}\frac{h_e}{\lambda_e} \Big\|\Big\{\frac{\partial^2(u-\Pi_h^{\nabla} I_h u)}{\partial \boldsymbol{n}_e^2}\Big\}\Big\|_e^2\Big)^{1/2} J_1(v_h, v_h)^{1/2}.
\end{align*}
The trace inequality, the inverse inequality on polynomials, the boundedness of $\Pi_h^{\nabla}$ in \eqref{Pi}, and $\lambda_e \eqsim 1$ give
\begin{align*}
& \Big(\sum_{e \in \mathcal{E}_h}\frac{h_e}{\lambda_e} \Big\|\Big\{\frac{\partial^2(u-\Pi_h^{\nabla} I_h u)}{\partial \bm{n}_e^2}\Big\}\Big\|_e^2\Big)^{1/2} \\
& \lesssim \sum_{e \in \mathcal{E}_h}h_e\Big(\Big\|\Big\{\frac{\partial^2 (u-u_\pi )}{\partial \bm{n}_e^2}\Big\}\Big\|_e^2+\Big\|\Big\{\frac{\partial^2 (u_\pi-\Pi_h^{\nabla} I_h u )}{\partial \bm{n}_e^2}\Big\}\Big\|_e^2\Big)^{1/2}\\
& \lesssim \sum_{K \in \mathcal{T}_h} ( |u-u_\pi |_{2, K}^2+h_K^2 |u-u_\pi |_{3, K}^2+ |\Pi_h^{\nabla} (I_h u-u_\pi ) |_{2, K}^2 )^{1/2}\\
& \lesssim |u -I_h u|_{2,h} + |u-u_\pi|_{2,h} + h |u-u_\pi|_{3,h}.
\end{align*}
This along with the Dupont-Scott theory and Lemma \ref{lem:interpIPVEM} yields
\begin{align*}
E_{A1}^2
& \lesssim \varepsilon^2  (|u -I_h u|_{2,h} + |u-u_\pi|_{2,h} + h |u-u_\pi|_{3,h}) J_1(v_h,v_h)^{1/2} \\
& \lesssim \varepsilon h^{s-2}|u|_s\| v_h \|_{\varepsilon,h}.
	\end{align*}
Collecting above estimates to derive
\begin{align*}
E_{A1} = E_{A1}^1 + E_{A1}^2 \lesssim  \varepsilon h^{s-2}|u|_s \| v_h \|_{\varepsilon,h}.
\end{align*}
	
	For $E_{A2}$, similar to the arguments in \cite[Lemma 5.3]{Zhao-Zhang-Chen-2018}, one has
	\[E_{A2} \lesssim \varepsilon^2 h^{s-2} |u|_s|v_h|_{2,h} \lesssim \varepsilon h^{s-2}|u|_s \| v_h \|_{\varepsilon,h}. \]

We next estimate
	\[E_{A3}  = -\varepsilon^2 \sum_{K\in \mathcal{T}_h}(\nabla \Delta u,\nabla v_h)_K +  \sum_{K\in \mathcal{T}_h}(\nabla u,  \nabla I_h^c v_h)_K - (f, v_h) ,	\]
  where $u$ is the exact solution of \eqref{originalSP}. For $k=2$, we have $ s = 3$. The integration by parts for $u \in H_0^2(\Omega)$ gives
\begin{align*}
0 & = \varepsilon^2 (\Delta^2 u, I_h^c v_h) -(\Delta u, I_h^c v_h)-(f,I_h^c v_h)\\
  & = -\varepsilon^2\sum_{K\in \mathcal{T}_h}(\nabla(\Delta u),\nabla I_h^c v_h )_K+\sum_{K\in \mathcal{T}_h}(\nabla u,\nabla I_h^c v_h )_K-(f,I_h^c v_h),
		\end{align*}
hence
\begin{equation}\label{eq:EA3}
E_{A3} =\sum_{K\in \mathcal{T}_h}\varepsilon^2\big(\nabla \Delta u,\nabla( I_h^c v_h-v_h)\big)_K+(f,I_h^c v_h-v_h).
\end{equation}
By the interpolation error estimate,
\[
E_{A3} \lesssim  (\varepsilon h|u|_3+h\|f\| )\| v_h \|_{\varepsilon,h}
 = h^{s-2} (\varepsilon  |u|_3 +\|f\|_{s-3} )\| v_h \|_{\varepsilon,h}.
\]
For $k\geq3$, by the integration by parts,
\begin{align*}
E_{A3}&= -\sum_{e\in \mathcal{E}_h}\varepsilon^2\int_e\frac{\partial \Delta u}{\partial \boldsymbol{n}_e}[v_h]\mathrm{d} s+\int_{e}\frac{\partial u}{\partial \bm n_e}[I_h^cv_h]\mathrm{d} s\notag\\
&= -\sum_{e\in \mathcal{E}_h}\varepsilon^2\int_e\Big(\frac{\partial \Delta u}{\partial \boldsymbol{n}_e}-\Pi_{0,e}^{k-3}\frac{\partial \Delta u}{\partial \boldsymbol{n}_e}\Big)[v_h-I_h^c v_h]\mathrm{d} s\\
&\lesssim \sum_{K\in \mathcal{T}_h} \varepsilon^2
h_K^{s-7/2}|u|_{s,K}h_K^{3/2}|v_h|_{2,K}\lesssim \varepsilon h^{s-2}|u|_{s}\| v_h \|_{\varepsilon,h}.
\end{align*}

Step 2: The remaining is to discuss $\varepsilon^2J_1(I_hu,v_h)$ and $\varepsilon^2J_3(I_hu,v_h)$. For $\varepsilon^2J_1(I_hu,v_h)$, the continuity of $u$ leads to
\begin{align*}
\varepsilon^2 J_1(I_h u,v_h)
& = \varepsilon^2 \sum_{e\in\mathcal{E}_h}\frac{\lambda_e}{|e|}\int_e\Big[\frac{\partial \Pi_h^{\nabla}I_hu }{\partial \bm n_e}\Big]
		\Big[\frac{\partial \Pi_h^{\nabla}v_h}{\partial \bm n_e}\Big] \mathrm{d} s \\
& = \varepsilon^2 \sum_{e\in\mathcal{E}_h}\frac{\lambda_e}{|e|}\int_e\Big[\frac{\partial (\Pi_h^{\nabla}I_hu-u)}{\partial \bm n_e}\Big]
		\Big[\frac{\partial \Pi_h^{\nabla}v_h}{\partial \bm n_e}\Big] \mathrm{d} s \\
& = \varepsilon^2 \sum_{e\in\mathcal{E}_h}\frac{\lambda_e^{1/2}}{|e|^{1/2}} \Big\| \Big[\frac{\partial (\Pi_h^{\nabla}I_hu-u)}{\partial \bm n_e}\Big] \Big\|_{0,e} J_1(v_h,v_h)^{1/2} .
\end{align*}
This further gives
\begin{align*}
\varepsilon^2 J_1(I_h u,v_h)
& \lesssim \varepsilon\sum_{K\in\mathcal{T}_h}(|\Pi_h^{\nabla}I_hu-u|_{2,K} + h_K^{-1}|\Pi_h^{\nabla}I_hu-u|_{1,K})
		\| v_h \|_{\varepsilon,h} \\
& \lesssim \varepsilon\sum_{K\in\mathcal{T}_h}\Big(|u - I_hu|_{2,K} + |u-u_\pi|_{2,K} + h_K^{-1}(|u-I_hu|_{1,K} + |u-u_\pi|_{1,K} )\Big)
		\| v_h \|_{\varepsilon,h}\\
& \lesssim \varepsilon h^{s-2}|u|_s \| v_h \|_{\varepsilon,h}.
	\end{align*}

For $\varepsilon^2J_3(I_hu,v_h)$, the continuity of $u$ and the trace inequality give
	\begin{align*}
& \varepsilon^2J_3(I_hu,v_h) \\
&  = \varepsilon^2\sum_{e\in\mathcal{T}_h}\int_e\Big[\frac{\partial (\Pi_h^{\nabla}I_hu-u)}{\partial \bm n_e}\Big]\Big[\frac{\partial^2\Pi_h^{\nabla}v_h}{\partial \bm n_e^2}\Big]\mathrm{d} s \lesssim \varepsilon^2\sum_{e\in \mathcal{T}_h}\Big\|\Big[\frac{\partial( \Pi_h^{\nabla}I_hu-u)}{\partial \bm n_e}\Big] \Big\|_{0,e}\Big\|\Big[\frac{\partial^2\Pi_h^{\nabla}v_h}{\partial \bm n_e^2}\Big] \Big\|_{0,e}\\
& \lesssim \varepsilon^2 \sum_{K\in\mathcal{T}_h}\Big(h_K^{1/2}|\Pi_h^{\nabla}I_hu-u|_{2,K} + h_K^{-1/2}|\Pi_h^{\nabla}I_hu-u|_{1,K}\Big) \Big(|\Pi_h^{\nabla}v_h|_{2,K}+|\Pi_h^{\nabla}v_h|_{2,K}^{1/2}|\Pi_h^{\nabla}v_h|_{3,K}^{1/2}\Big)\\
&\lesssim \varepsilon \sum_{K\in \mathcal{T}_h}\Big( |u-I_hu|_{2,K}+|u-u_{\pi}|_{2,k}+h_K^{-1}(|u-I_hu|_{1,k}+|u-u_{\pi}|_{1,K})\Big)\| v_h \|_{\varepsilon,h}\\
& \lesssim \varepsilon h^{s-2}|u|_s \| v_h \|_{\varepsilon,h}.
	\end{align*}

Combining the estimates in the above two steps, we immediately obtain
\[E_{A}(u,v_h)+E_J(u,v_h)\lesssim h^{s-2} (\varepsilon |u|_s +  \|f\|_{s-3} ) \| v_h \|_{\varepsilon,h},\]
which completes the proof.
\end{proof}

To sum up the above results, we obtain the error estimate for the IPVEM described as follows.
\begin{theorem}\label{estimate}
Given $k \ge 2$ and $f \in H^{s-3}(\Omega)$ with $3 \leq s \leq$ $k+1$,
assume that $u \in H_0^2(\Omega) \cap H^s(\Omega)$ is the exact solution of \eqref{originalSP}. Then there holds
\begin{equation}
	\|u-u_h\|_{\varepsilon,h}\lesssim \left\{\begin{array}{l}
		 ( h^{s-1} + \varepsilon h^{s-2}) |u|_s + h^{s-2}\|f\|_{s-3}, \\
	     ( h^{s-1} + \varepsilon h^{s-2}) |u|_s + h^{s-2}\|f\|_{s-3}.
	\end{array}\right.
\end{equation}
\end{theorem}
\begin{proof}
We only need to bound each term of the right-hand side of \eqref{Strangerr} in Lemma \ref{lem:StrangIPVEM}.  By by the trace inequality, the boundedness of $\Pi_h^{\nabla}$ in \eqref{Pi}, the interpolation error estimate in Lemma \ref{lem:interpIPVEM} and the Dupont-Scott theory, there exists $u_\pi\in \mathbb{P}_k(\mathcal{T}_h)$ such that
\[\| u - I_h u \|_{\varepsilon,h} + \| u - u_\pi \|_{\varepsilon,h} \lesssim h^{\ell-1}|u|_\ell + \varepsilon h^{\ell'-2}|u|_{\ell'}, \qquad 2\le \ell, \ell' \le k+1.\]
If $\ell = \ell' = s$, then
\[\| u - I_h u \|_{\varepsilon,h} + \| u - u_\pi \|_{\varepsilon,h} \lesssim ( h^{s-1} + \varepsilon h^{s-2}) |u|_s.\]
If $\ell = s-1$, $\ell' = s$, then
\[\| u - I_h u \|_{\varepsilon,h} + \| u - u_\pi \|_{\varepsilon,h} \lesssim h^{s-2} ( |u|_{s-1} + \varepsilon |u|_s).\]
The proof is completed by combining the above equations, \eqref{errRhs} and Lemma \ref{errConsistency}.
\end{proof}

\subsection{Uniform error estimate in the lowest order case}

Let $u^0$ be the solution of the following boundary value problem:
\begin{equation}\label{bdp}
	\left\{\begin{array}{ll}
		-\Delta u^0 = f & \quad {\rm in} \quad \Omega\\
		u^0=0 &  \quad {\rm on} \quad \partial \Omega.
	\end{array}\right.
\end{equation}
The following regularity is well-known and can be found in \cite{NTW01} for instance.
\begin{lemma}\label{lem:regular}
	If $\Omega$ is a bounded convex polygonal domain, then
	\[
	|u|_2+\varepsilon|u|_3\lesssim \varepsilon^{-1/2}\|f\|\qquad \text{and}\qquad |u-u^0|_1\lesssim \varepsilon^{1/2}\|f\|
	\]
	for all $f\in L^2(\Omega)$.
\end{lemma}

\begin{theorem}
Let $k=2$ be the order of the virtual element space.
Assume that $f \in L^2(\Omega)$ and $u \in H_0^2(\Omega)\cap H^3(\Omega)$ is the exact solution to \eqref{originalSP}.
If $\Omega$ is a bounded convex polygonal domain, then
\begin{equation}\label{result:consistent}
	\|u-u_h\|_{\varepsilon,h} \lesssim h^{1/2} \|f\| .
\end{equation}
\end{theorem}
\begin{proof}
We only need to estimate the right-hand side of \eqref{Strangerr} in Lemma \ref{lem:StrangIPVEM} term by term.

Step 1:
By the multiplicative trace inequality, the boundedness of $\Pi_h^{\nabla}$
and $\lambda_e \eqsim 1$ and Lemma \ref{lem:regular}, one gets
\begin{align}
\varepsilon^2J_1(u-I_hu,u-I_hu)\lesssim h\|f\|^2,\quad
\varepsilon^2J_1(u-u_{\pi},u-u_{\pi})\lesssim h\|f\|^2.\label{j1}
\end{align}
Following the same manipulations in Theorem 4.3 of \cite{ZZC20}, combing with \eqref{j1}, one may obtain
\begin{align}
\|u-I_hu\|_{\varepsilon,h} + \|u-u_{\pi}\|_{\varepsilon,h} \lesssim h^{1/2}\|f\| \label{half12}.
\end{align}
The estimate of the right-hand side is given in \eqref{errRhs} with $s = 3$ in this theorem.

Step 2: For the consistency term, we first consider $E_{A1}, E_{A2}, E_{A3}$ defined in Lemma \ref{rem:Strangdecomp}.
By the trace inequalities in Lemma \ref{lem:trace},
 \begin{align}
E_{A1}^1
& = \varepsilon^2\sum_{e \in \mathcal{E}_h}\int_e  \Big\{\frac{\partial^2 (u-\Pi_h ^\nabla I_h u)}{\partial \boldsymbol{n}_e^2}\Big\} \Big[\frac{\partial (v_h - \Pi_h^\nabla v_h)}{\partial \boldsymbol{n}_e}\Big]  \d s \nonumber\\
& \lesssim \varepsilon^2\sum_{e \in \mathcal{E}_h}\Big\| \Big\{\frac{\partial^2 (u-\Pi_h ^\nabla I_h u)}{\partial \boldsymbol{n}_e^2}\Big\}\Big\|_{0,e}\Big\|\Big[\frac{\partial (v_h - \Pi_h^\nabla v_h)}{\partial \boldsymbol{n}_e}\Big]  \Big\|_{0,e} \nonumber\\
& \lesssim \varepsilon^2\sum_{K \in \mathcal{T}_h} \|u-\Pi_h ^\nabla I_h u\|_{2,K}^{1/2}\|u-\Pi_h ^\nabla I_h u\|_{3,K}^{1/2} \nonumber \\
& \qquad \times \Big( h_K^{1/2}|v_h - \Pi_h^\nabla v_h|_{2,K} + h_K^{-1/2} |v_h - \Pi_h^\nabla v_h|_{1,K}\Big). \label{EA11uniform}
 \end{align}
This along with the boundedness of $\Pi_h^{\nabla}$, the interpolation error estimate and Lemma \ref{lem:regular} gives
  \begin{align*}
E_{A1}^1
& \lesssim h^{1/2} \varepsilon^2  |u|_2^{1/2} |u|_3^{1/2} |v_h|_{2,h} \lesssim h^{1/2}\|f\|\| v_h \|_{\varepsilon,h}.
 \end{align*}
By the boundedness of $\Pi_h^{\nabla}$, we similarly obtain
\[
E_{A1}^2 = \varepsilon^2\sum_{e \in \mathcal{E}_h}\int_e  \Big\{\frac{\partial^2 (u-\Pi_h ^\nabla I_h u)}{\partial \boldsymbol{n}_e^2}\Big\} \Big[\frac{\partial \Pi_h^\nabla v_h}{\partial \boldsymbol{n}_e}\Big]  \d s\lesssim h^{1/2}\|f\|\| v_h \|_{\varepsilon,h}.
\]
Therefore,
\[E_{A1} = E_{A1}^1 + E_{A1}^2 \lesssim h^{1/2}\|f\|\| v_h \|_{\varepsilon,h}.\]

For $E_{A2}$, we first obtain from the continuity of $v_h$ at vertices that
\begin{equation*}\label{exactquad}
\int_{e}q\Big[\frac{\partial v_h}{\partial \bm t_e}\Big] \d s=-\int_e\frac{\partial q}{\partial \bm t_e}[v_h]=0, \quad q\in \mathbb{P}_{k-2}(e),\quad e\in \mathcal{E}_h.
\end{equation*}
By setting $w|_e=\dfrac{\partial^2 u}{\partial \bm{n}_e \partial \bm{t}_e}$, where $w$ is extended outside $e$ so that it is constant along the lines perpendicular to $e$, this implies
 \begin{align*}
E_{A2}
& =\varepsilon^2\sum_{e \in \mathcal{E}_h}\int_e \Big(\frac{\partial^2 u}{\partial \boldsymbol{n}_e \partial \boldsymbol{t}_e} - \Pi_{0,e}^0\Big(\frac{\partial^2 u}{\partial \boldsymbol{n}_e \partial \boldsymbol{t}_e}\Big)\Big)\Big(\left[\frac{\partial v_h}{\partial \boldsymbol{t}_e}\right]-\Pi_{0,e}^0\left[\frac{\partial v_h}{\partial \boldsymbol{t}_e}\right]\Big) \d s\\
& \le \varepsilon^2\sum_{e \in \mathcal{E}_h}\| w - \Pi_{0,e}^0w \|_{0,e} \|[\partial_{{\bm t}_e} v_h] -\Pi_{0,e}^0[\partial_{{\bm t}_e} v_h]\|_{0,e} \\
& \le \varepsilon^2\sum_{e \in \mathcal{E}_h}\| w - \Pi_{0,K}^0w \|_{0,e} \|[\partial_{{\bm t}_e} v_h] -\Pi_{0,e}^0[\partial_{{\bm t}_e} v_h]\|_{0,e},
 \end{align*}
where we have used the minimization property of $L^2$ projection and $e$ is an edge of $K$. As done in the last step of \eqref{EA11uniform}, we immediately obtain
\[E_{A2}\lesssim  h^{1/2}\|f\|\| v_h \|_{\varepsilon,h}.\]

For $E_{A3}$, we choose $k=2$ with the formula given in \eqref{eq:EA3}, i.e.,
\[E_{A3}=\sum_{K\in \mathcal{T}_h}\varepsilon^2\big(\nabla \Delta u,\nabla( I_h^c v_h-v_h)\big)_K+(f,I_h^c v_h-v_h),\]
where
\[|(f,I_h^c v_h-v_h)| \lesssim h \|f\| \| v_h \|_{\varepsilon,h}\]
is the direct consequence of the interpolation error estimate.
We obtain from the Cauchy-Schwarz inequality, the interpolation error estimate and the boundeness of $\Pi_h^{\nabla}$ that
 \begin{align*}
& \Big|\sum_{K\in \mathcal{T}_h}\big(\nabla \Delta u,\nabla( I_h^c v_h-v_h)\big)_K\Big| \\
& = \Big| \sum_{K\in\mathcal{T}_h}\int_{\partial K}\Delta u\frac{\partial(I_h^c v_h-v_h)}{\partial \bm n}\mathrm{d} s-\int_K \Delta u \Delta (I_h^c v_h-v_h)\mathrm{d} x  \Big|\\
& \le \sum_{K\in \mathcal{T}_h}\|\Delta u\|_{0,\partial K}\Big\|\frac{\partial(I_h^c v_h-v_h)}{\partial \bm n_e}\Big\|_{0,\partial K}+\sum_{K\in \mathcal{T}_h}\|\Delta u\|_{0,K}|I_h^c v_h-v_h|_{2,K}\\
& \lesssim \Big(h^{1/2}\|\Delta u\|_{0,\Omega}^{1/2}\|\Delta u\|_{1,\Omega}^{1/2}+\|\Delta u\|_{0,\Omega}\Big)|v_h|_{2,h}.
 \end{align*}
This implies
\[
\varepsilon^2 \Big|\sum_{K\in \mathcal{T}_h}\big(\nabla \Delta u,\nabla( I_h^c v_h-v_h)\big)_K\Big|
\lesssim ( h^{1/2} + \varepsilon^{1/2}) \|f\| \|v_h\|_{\varepsilon,h}
\lesssim h^{1/2} \|f\| \|v_h\|_{\varepsilon,h}
\]
if $\varepsilon \le h$.  On the other hand, for $\varepsilon \ge h$,
 \begin{align*}
 	\varepsilon^2 \Big|\sum_{K\in \mathcal{T}_h}\big(\nabla \Delta u,\nabla( I_h^c v_h-v_h)\big)_K\Big|\lesssim &\varepsilon h|u|_3 \|v_h\|_{\varepsilon,h}
 	\lesssim  h\varepsilon^{-1/2}\|f\|\|v_h\|_{\varepsilon,h}
 	\lesssim  h^{1/2}\|f\|\|v_h\|_{\varepsilon,h} .
 \end{align*}

Step 3: The remaining consistency terms of $\varepsilon^2J_1(I_hu,v_h)$ and $\varepsilon^2J_3(I_hu,v_h)$ can be bounded as in the second step, so we omit the details, with the result given by
\[\varepsilon^2 ( J_1(I_hu,v_h) + J_3(I_hu,v_h)) \lesssim  h^{1/2} \|f\| \|v_h\|_{\varepsilon,h}.\]

The uniform estimate follows by combining the above bounds.	
\end{proof}

\section{Numerical examples} \label{sec:numerical}

In this section, we report the performance with several examples by testing the accuracy and the robustness with respect to the singular parameter $\varepsilon$.
For simplicity, we only consider the lowest-order element ($k = 2$) and the domain $\Omega$ is taken as the unit square $(0,1)^2$.

\begin{example}\label{example1}
The source term is chosen in such a way that the exact solution is
  \[ u = 10 x^2 y^2 (1-x)^2 (1-y)^2 \sin(\pi x). \]
\end{example}

Let $u$ be the exact solution of \eqref{originalSP} and $u_h$ the discrete solution of the underlying VEM \eqref{IPVEM1}. Since the VEM solution is not explicitly known inside the polygonal elements, we will evaluate the errors by comparing the exact solution $u$ with the interpolation $I_h u_h$. In this way, the discrete error in terms of the discrete energy norm is quantified by
\[E_I = \Big( \sum\limits_{K \in \mathcal{T}_h} (\varepsilon^2 |u -  I_h u_h|_{2,K}^2 + |u - I_h u_h|_{1,K}^2 ) \Big)^{1/2}.\]

To test the accuracy of the proposed method we consider a sequence of meshes, which is a Centroidal Voronoi Tessellation of the unit square in 32, 64, 128, 256 and 512 polygons. These meshes are generated by the MATLAB toolbox - PolyMesher introduced in \cite{Talischi-Paulino-Pereira-2012}.
 The convergence orders of the errors against the mesh size $h$ are shown in Table \ref{tab:exam1poly}.
As observed from Table \ref{tab:exam1poly}, the VEM ensures the linear convergence as $\varepsilon\to 0$, which is consistent with the theoretical prediction in Theorem \ref{estimate}. Moreover, a stable trend of the errors is observed as $\varepsilon$ decreases to zero.

\begin{table}[!htb]
  \centering
  \caption{The convergence rate for Example \ref{example1}}\label{tab:exam1poly}
  \begin{tabular}{ccccccccccccccccc}
  \toprule[0.2mm]
  $\varepsilon \backslash N$ & 32   &64   &128   &256   &512  & Rate\\
  \midrule[0.3mm]
   1e-0  & 7.7401e-01   & 4.9029e-01   & 2.7904e-01   & 1.6441e-01   & 1.1698e-01 & 1.41\\
   1e-1  & 7.2306e-02   & 4.8619e-02   & 2.8865e-02   & 1.7533e-02   & 1.2075e-02 & 1.33\\
   1e-2  & 2.4167e-02   & 1.8044e-02   & 1.3098e-02   & 9.4385e-03   & 6.6785e-03 & 0.93\\
   1e-3  & 2.3908e-02   & 1.7905e-02   & 1.3125e-02   & 9.4839e-03   & 6.7673e-03 & 0.91\\
   1e-4  & 2.3910e-02   & 1.7910e-02   & 1.3132e-02   & 9.4912e-03   & 6.7754e-03 & 0.91\\
   1e-5  & 2.3910e-02   & 1.7910e-02   & 1.3133e-02   & 9.4912e-03   & 6.7755e-03 & 0.91\\
  \bottomrule[0.2mm]
\end{tabular}
\end{table}

\begin{figure}[!htb]
  \centering
  \includegraphics[scale=0.7,trim = 50 0 50 0,clip]{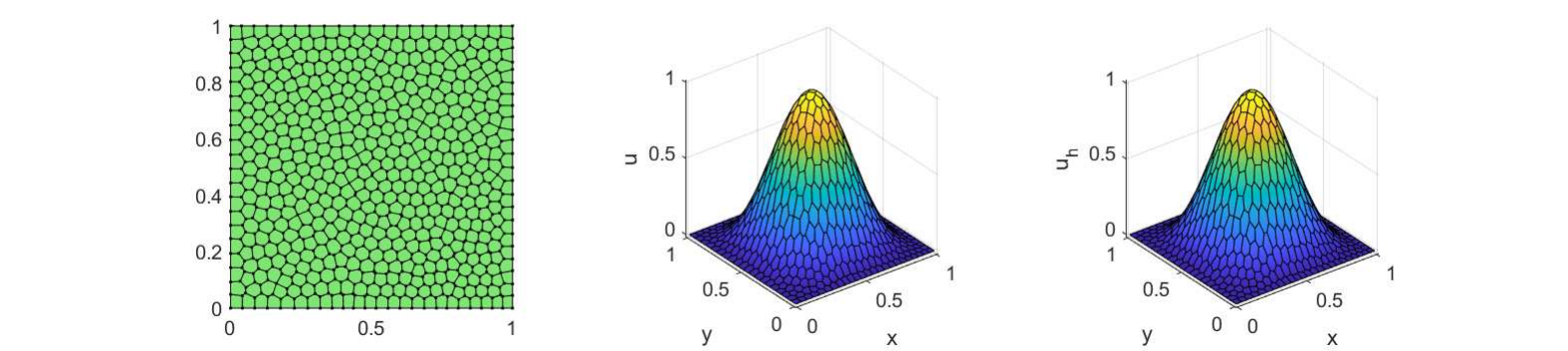}\\
  \caption{Numerical and exact solutions for Example \ref{example2}~($\varepsilon = 10^{-10}$)}\label{fig:example2}
\end{figure}

\begin{example}\label{example2}
The exact solution is given by $u = \sin(\pi x)^2\sin(\pi y)^2$.
\end{example}

\begin{figure}[!htb]
  \centering
  \includegraphics[scale=0.5]{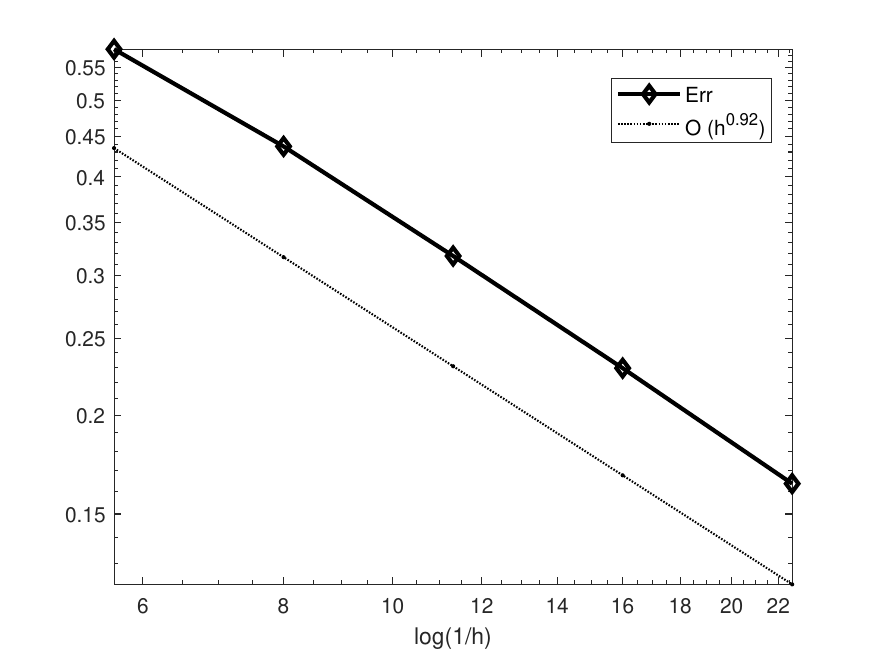}\\
  \caption{Convergence rate for Example \ref{example2}~($\varepsilon = 10^{-10}$)}\label{fig:example2rate}
\end{figure}

Fig.~\ref{fig:example2} displays the numerical solution and exact solution for $\varepsilon = 10^{-10}$. These two figures fit well, and the optimal first-order convergence is still observed in Fig.~\ref{fig:example2rate} although the singular parameter $\varepsilon$ is very small.

\section{Conclusion}

We proposed a modified interior penalty virtual element method to solve the fourth-order singular perturbation problem. Drawing inspiration from the technique of the modified Morley element, we deduced optimal convergence in the energy norm and demonstrated robustness with the perturbation parameter in the lowest-order case.

\section*{CRediT authorship contribution statement}

Fang Feng and Yue Yu collaborated closely to shape the conceptualization, methodology, and writing of this research. Additionally, Yue Yu took charge of implementing the discrete method in the study.

\section*{Declaration of competing interest}

The authors declare that they have no known competing financial interests or personal relationships that could have appeared to
influence the work reported in this paper.

\section*{Data availability}

Data will be made available on request.

\section*{Acknowledgements}

Yue Yu was partially supported by the National Science Foundation for Young Scientists of China (No. 12301561).

\bibliographystyle{plain} 
\bibliography{Refs_IPVEM}

\begin{thebibliography}{10}

\bibitem{Adak-Nataraj-2020}
D.~Adak and S.~Natarajan.
\newblock Virtual element method for semilinear sine-{G}ordon equation over
  polygonal mesh using product approximation technique.
\newblock {\em Math. Comput. Simulation}, 172:224--243, 2020.

\bibitem{Ahmad-Alsaedi-Brezzi-2013}
B.~Ahmad, A.~Alsaedi, F.~Brezzi, L.D. Marini, and A.~Russo.
\newblock Equivalent projectors for virtual element methods.
\newblock {\em Comput. Math. Appl.}, 66(3):376--391, 2013.

\bibitem{Alvarez-Beirao-Dassi-2023}
S.~N. Alvarez, L.~{Beir{\~{a}}o Da Veiga}, F.~Dassi, V.~Gyrya, and G.~Manzini.
\newblock The virtual element method for a {2D} incompressible {MHD} system.
\newblock {\em Math. Comput. Simulation}, 211:301--328, 2023.

\bibitem{Antonietti-Beirao-Manzini2022}
P.~F. Antonietti, L.~Beir\~{a}o~da Veiga, and G.~Manzini.
\newblock {\em The virtual element method and its applications}.
\newblock Springer, Cham, 2022.

\bibitem{Antonietti-Bruggi-Scacchi-17}
P.~F. Antonietti, M.~Bruggi, S.~Scacchi, and M.~Verani.
\newblock On the virtual element method for topology optimization on polygonal
  meshes: a numerical study.
\newblock {\em Comput. Math. Appl.}, 74(5):1091--1109, 2017.

\bibitem{Antonietti-Manzini-Verani-2018}
P.~F. Antonietti, G.~Manzini, and M.~Verani.
\newblock The fully nonconforming virtual element method for biharmonic
  problems.
\newblock {\em Math. Models Methods Appl. Sci.}, 28(2):387--407, 2018.

\bibitem{Beirao-Brezzi-Cangiani-2013}
L.~{Beir{\~{a}}o Da Veiga}, F.~Brezzi, A.~Cangiani, G.~Manzini, L.~D. Marini,
  and A.~Russo.
\newblock Basic principles of virtual element methods.
\newblock {\em Math. Models Meth. Appl. Sci.}, 23(1):199--214, 2013.

\bibitem{Beirao-Brezzi-Marini-2014}
L.~{Beir{\~{a}}o Da Veiga}, F.~Brezzi, L.~D. Marini, and A.~Russo.
\newblock The {H}itchhiker's guide to the virtual element method.
\newblock {\em Math. Models Meth. Appl. Sci.}, 24(8):1541--1573, 2014.

\bibitem{Beirao-Dassi-Manzini-2023}
L.~{Beir{\~{a}}o Da Veiga}, F.~Dassi, G.~Manzini, and L.~Mascotto.
\newblock The virtual element method for the {3D} resistive magnetohydrodynamic
  model.
\newblock {\em Math. Models Methods Appl. Sci.}, 33(3):643--686, 2023.

\bibitem{Beirao-Lovadina-Vacca-2018}
L.~{Beir{\~{a}}o Da Veiga}, C.~Lovadina, and G.~Vacca.
\newblock Virtual elements for the {N}avier-{S}tokes problem on polygonal
  meshes.
\newblock {\em SIAM J. Numer. Anal.}, 56(3):1210--1242, 2018.

\bibitem{Beirao-Mora-Vacca-2019}
L.~{Beir{\~{a}}o Da Veiga}, D.~Mora, and G.~Vacca.
\newblock The {S}tokes complex for virtual elements with application to
  {N}avier–{S}tokes flows.
\newblock {\em J. Sci. Comput.}, 81:990--1018, 2019.

\bibitem{Brenner2011NeilanC0IP}
S.~C. Brenner and M.~Neilan.
\newblock A {$C^0$} interior penalty method for a fourth order elliptic
  singular perturbation problem.
\newblock {\em SIAM J. Numer. Anal.}, 49:869--892, 2011.

\bibitem{BS2008}
S.~C. Brenner and L.~R. Scott.
\newblock {\em The Mathematical Theory of Finite Element Methods}.
\newblock Springer, New York, 2008.

\bibitem{Brenner2005SungC0IP}
S.~C. Brenner and L.~Sung.
\newblock {$C^0$} interior penalty methods for fourth order elliptic boundary
  value problems on polygonal domains.
\newblock {\em J. Sci. Comput.}, 22/23:83--118, 2005.

\bibitem{Brezzi-Buffa-Lipnikov-2009}
F.~Brezzi, A.~Buffa, and K.~Lipnikov.
\newblock Mimetic finite differences for elliptic problems.
\newblock {\em M2AN Math. Model. Numer. Anal.}, 43(2):277--295, 2009.

\bibitem{Brezzi-Marini-2013}
F.~Brezzi and L.~D. Marini.
\newblock Virtual element methods for plate bending problems.
\newblock {\em Comput. Methods Appl. Mech. Engrg.}, 253:455--462, 2013.

\bibitem{Caceres-Gatica-2017}
E.~C\'{a}ceres and G.~N. Gatica.
\newblock A mixed virtual element method for the pseudostress-velocity
  formulation of the {S}tokes problem.
\newblock {\em IMA J. Numer. Anal.}, 37:296--331, 2017.

\bibitem{Cangiani-Manzini-Sutton-2016}
A.~Cangiani, G.~Manzini, and O.~J. Sutton.
\newblock Conforming and nonconforming virtual element methods for elliptic
  problems.
\newblock {\em IMA J. Numer. Anal.}, 37(3):1317--1354, 2016.

\bibitem{Chen-HuangJ-2018}
L.~Chen and J.~Huang.
\newblock Some error analysis on virtual element methods.
\newblock {\em Calcolo}, 55(1):5, 2018.

\bibitem{Chen-HuangX-2020}
L.~Chen and X.~Huang.
\newblock {Nonconforming virtual element method for {$2m$}-th order partial
  differential equations in {$R^n$}}.
\newblock {\em Math. Comput.}, 89(324):1711--1744, 2020.

\bibitem{Chi-Pereira-Menezes-Paulino-20}
H.~Chi, A.~Pereira, I.~F.~M. Menezes, and G.~H. Paulino.
\newblock Virtual element method ({VEM})-based topology optimization: an
  integrated framework.
\newblock {\em Struct. Multidiscip. Optim.}, 62(3):1089--1114, 2020.

\bibitem{Chinosi-Marini-2016}
C.~Chinosi and L.~D. Marini.
\newblock Virtual element method for fourth order problems: {$L^2$}-estimates.
\newblock {\em Comput. Math. Appl.}, 72(8):1959--1967, 2016.

\bibitem{DeDios-Lipnikov-Manzini-2016}
B.~A. De~Dios, K.~Lipnikov, and G.~Manzini.
\newblock The nonconforming virtual element method.
\newblock {\em ESAIM Math. Model. Numer. Anal.}, 50(3):879--904, 2016.

\bibitem{FHH19}
F.~Feng, W.~Han, and J.~Huang.
\newblock Virtual element method for an elliptic hemivariational inequality
  with applications to contact mechanics.
\newblock {\em J. Sci. Comput.}, 81(4036991):2388--2412, 2019.

\bibitem{FHH22}
F.~Feng, W.~Han, and J.~Huang.
\newblock A nonconforming virtual element method for a fourth-order
  hemivariational inequality in kirchhoff plate problem.
\newblock {\em J. Sci. Comput.}, 90(3):89, 24 pp., 2022.

\bibitem{Gatica-Munar-2018}
G.~N. Gatica and M.~Munar.
\newblock A mixed virtual element method for the {N}avier–{S}tokes equations.
\newblock {\em Math. Models Methods Appl. Sci.}, 28(14):2719–2762, 2018.

\bibitem{Huang2021YuMedius}
J.~Huang and Y.~Yu.
\newblock A medius error analysis for nonconforming virtual element methods for
  poisson and biharmonic equations.
\newblock {\em J. Comput. Appl. Math.}, 386:113229, 20 pp, 2021.

\bibitem{Ling-wang-han-2020}
M.~Ling, F.~Wang, and W.~Han.
\newblock The nonconforming virtual element method for a stationary {S}tokes
  hemivariational inequality with slip boundary condition.
\newblock {\em J. Sci. Comput.}, 85(3):Paper No. 56, 19, 2020.

\bibitem{Liu-Li-Chen-2017}
X.~Liu, J.~Li, and Z.~Chen.
\newblock A nonconforming virtual element method for the {S}tokes problem on
  general meshes.
\newblock {\em Comput. Methods Appl. Mech. Engrg}, 320:694--711, 2017.

\bibitem{NTW01}
T.~K. Nilssen, X.~Tai, and R.~Winther.
\newblock A robust nonconforming {$H^2$}-element.
\newblock {\em Math. Comput.}, 70(234):489--505, 2001.

\bibitem{Carsten2023IP}
B.~Philipp, C.~Carsten, and S.~Julian.
\newblock Local parameter selection in the {$C^0$} interior penalty method for
  the biharmonic equation.
\newblock {\em arXiv:2209.05221v2}, 2023.

\bibitem{qiu-wang-ling-zhao-2023}
J.~Qiu, F.~Wang, M.~Ling, and J.~Zhao.
\newblock The interior penalty virtual element method for the fourth-order
  elliptic hemivariational inequality.
\newblock {\em Commun. Nonlinear Sci. Numer. Simul.}, 127(4644807):Paper No.
  107547, 17, 2023.

\bibitem{Semper1992}
B.~Semper.
\newblock Conforming finite element approximations for a fourth-order singular
  perturbation problem.
\newblock {\em SIAM J. Numer. Anal.}, 29(4):1043--1058, 1992.

\bibitem{Talischi-Paulino-Pereira-2012}
C.~Talischi, G.~H. Paulino, A.~Pereira, and Ivan F.~M. Menezes.
\newblock Polymesher: a general-purpose mesh generator for polygonal elements
  written in matlab.
\newblock {\em Struct. Multidiscip. Optim.}, 45(3):309--328, 2012.

\bibitem{Wang-Wei-2017}
F.~Wang and H.~Wei.
\newblock Virtual element method for simplified friction problem.
\newblock {\em Appl. Math. Lett.}, 85(3820290):125--131, 2018.

\bibitem{Wang-wu-han-2021}
F.~Wang, B.~Wu, and W.~Han.
\newblock The virtual element method for general elliptic hemivariational
  inequalities.
\newblock {\em J. Comput. Appl. Math.}, 389(4194398):Paper No. 113330, 19,
  2021.

\bibitem{wang-zhao-2021}
F.~Wang and J.~Zhao.
\newblock Conforming and nonconforming virtual element methods for a
  {K}irchhoff plate contact problem.
\newblock {\em IMA J. Numer. Anal.}, 41(2):1496--1521, 2021.

\bibitem{Wang01}
M.~Wang.
\newblock On the necessity and sufficiency of the patch test for convergence of
  nonconforming finite elements.
\newblock {\em SIAM J. Numer. Anal.}, 39(2):363--384, 2001.

\bibitem{WXH06}
M.~Wang, J.~Xu, and Y.~Hu.
\newblock Modified {M}orley element method for a fourth order elliptic singular
  perturbation problem.
\newblock {\em J. Comput. Math.}, 24(2):113--120, 2006.

\bibitem{xiao-ling-2023}
W.~Xiao and M.~Ling.
\newblock Virtual element method for a history-dependent
  variational-hemivariational inequality in contact problems.
\newblock {\em J. Sci. Comput.}, 96(3):Paper No. 82, 21, 2023.

\bibitem{Zhang-chi-Paulino-20}
X.~Zhang, H.~Chi, and G.~H. Paulino.
\newblock Adaptive multi-material topology optimization with hyperelastic
  materials under large deformations: a virtual element approach.
\newblock {\em Comput. Methods Appl. Mech. Engrg.}, 370(4129484):112976, 34,
  2020.

\bibitem{ZMZW2023IPVEM}
J.~Zhao, S.~Mao, B.~Zhang, and F.~Wang.
\newblock The interior penalty virtual element method for the biharmonic
  problem.
\newblock {\em Math. Comp.}, 92(342):1543--1574, 2023.

\bibitem{ZZC20}
J.~Zhao, B.~Zhang, and S.~Chen.
\newblock The nonconforming virtual element method for fourth-order singular
  perturbation problem.
\newblock {\em Adv. Comput. Math.}, 46(2):Paper No. 19, 23, 2020.

\bibitem{Zhao-Zhang-Chen-2018}
J.~Zhao, B.~Zhang, S.~Chen, and S.~Mao.
\newblock The {Morley-type} virtual element for plate bending problems.
\newblock {\em J. Sci. Comput.}, 76(1):610--629, 2018.

\end{thebibliography}

\end{document}